\documentclass[12pt,english,a4paper]{amsart}
\usepackage[T1]{fontenc}
\usepackage[latin1]{inputenc}
\usepackage{geometry}
\usepackage{amssymb}

\makeatletter
 \theoremstyle{plain}    
 \newtheorem{thm}{Theorem}[section]
 \numberwithin{equation}{section} 
 \numberwithin{figure}{section} 
 \theoremstyle{plain}
 \theoremstyle{remark}
 \newtheorem{rem}[thm]{Remark}
 \theoremstyle{remark}    
 \newtheorem{notation}[thm]{Notation} 
 \theoremstyle{plain}    
 \newtheorem{prop}[thm]{Proposition} 
 \theoremstyle{plain}    
 \newtheorem{lem}[thm]{Lemma} 
 \theoremstyle{plain}    
 \newtheorem{cor}[thm]{Corollary} 
 \theoremstyle{remark}    
 \newtheorem{acknowledgement}[thm]{Acknowledgement} 

\def\makebbb#1{
    \expandafter\gdef\csname#1\endcsname{
        \ensuremath{\Bbb{#1}}}
}
\makebbb{R}
\makebbb{N}
\makebbb{Z}
\makebbb{C}
\makebbb{G}
\makebbb{E}
\makebbb{H}
\makebbb{P}
\makebbb{K}
\makebbb{B}

\usepackage{babel}
\makeatother
\begin{document}

\title{Holomorphic Morse Inequalities on Manifolds with Boundary}

\author{Robert Berman}

\address{Department of Mathematics, Chalmers University of Technology, Eklandag.
86, SE-412 96 Göteborg}

\email{robertb@math.chalmers.se}

\keywords{Line bundles, Cohomology, Harmonic forms, Holomorphic sections, Bergman
kernel. \emph{MSC (2000):} 32A25, 32L10, 32L20}

\begin{abstract}
Let $X$ be a compact complex manifold with boundary and let $L^{k}$
be a high power of a hermitian holomorphic line bundle over $X.$
When $X$ has no boundary, Demailly's holomorphic Morse inequalities
give asymptotic bounds on the dimensions of the Dolbeault cohomology
groups with values in $L^{k},$ in terms of the curvature of $L.$
We extend Demailly's inequalities to the case when $X$ has a boundary
by adding a boundary term expressed as a certain average of the curvature
of the line bundle and the Levi curvature of the boundary. Examples
are given that show that the inequalities are sharp.
\end{abstract}
\maketitle

\section{Introduction}

Let $X$ be a compact $n-$dimensional complex manifold with boundary.
Let $\rho$ be a defining function of the boundary of $X,$ i.e. $\rho$
is defined in a neighborhood of the boundary of $X,$ vanishing on
the boundary and negative on $X.$ We take a hermitian metric $\omega$
on $X$ such that $d\rho$ is of unit-norm close to the boundary of
$X.$ The restriction of the two-form $i\partial\overline{\partial}\rho$
to the maximal complex subbundle $T^{1,0}(\partial X)$ of the tangent
bundle of $\partial X,$ is the Levi curvature form of the boundary
$\partial X.$ It will be denoted by $\mathcal{L}.$ Furthermore,
let $L$ be a hermitian holomorphic line bundle over $X$ with fiber
metric $\phi,$ so that $i\partial\overline{\partial}\phi$ is the
curvature two-form of $L.$ It will be denoted by $\Theta.$ The line
bundle $L$ is assumed to be smooth up to the boundary of $X.$ Strictly
speaking, $\phi$ is a collection of local functions. Namely, let
$s_{i}$ be a local holomorphic trivializing section of $L,$ then
locally, $\left|s_{i}(z)\right|^{2}=e^{-\phi_{i}(z)}.$ The notation
$\eta_{p}:=\eta^{p}/p!$ will be used in the sequel, so that the volume
form on $X$ may be written as $\omega_{n}.$ 

When $X$ is a compact manifold without boundary Demailly's (weak)
holomorphic Morse inequalities \cite{d1} give asymptotic bounds on
the dimension of the Dolbeault cohomology groups associated to the
$k:$th tensor power of the line bundle $L:$ \begin{equation}
\dim_{\C}H^{0,q}(X,L^{k})\leq k^{n}(-1)^{q}(\frac{1}{2\pi})^{n}\int_{X(q)}\Theta_{n}+o(k^{n}),\label{demailly weak ineaq}\end{equation}
 where $X(q)$ is the subset of $X$ where the curvature-two form
$\Theta$ has exactly $q$ negative eigenvalues, i.e the set where
$\textrm{index$($}\Theta)=q.$ Demailly's inspiration came from Witten's
analytical proof of the classical Morse inequalities for the Betti
numbers of a \emph{}real \emph{}manifold \cite{wi}, where the role
of the fiber metric $\phi$ is played by a Morse function. Subsequently,
holomorphic Morse inequalities on manifolds with boundary where studied.
The cases of $q-$convex and $q-$concave boundary were studied by
Bouche \cite{bo}, and Marinescu \cite{m1}, respectively, and they
obtained the same curvature integral as in the case when $X$ has
no boundary. However, it was assumed that, close to the boundary,
the curvature of the line bundle $L$ is adapted to the curvature
of the boundary. For example, on a pseudoconcave manifold (i.e the
Levi form is negative on the boundary) it is assumed that the curvature
of $L$ is non-positive close to the boundary. This is related to
the well-known fact that in the global $L^{2}-$ estimates for the
$\overline{\partial}-$operator of Morrey-Kohn-Hörmander-Kodaira there
is a curvature term from the line bundle as well from the boundary
and, in general, it is difficult to control the sign of the total
curvature contribution. 

In the present paper we will consider an arbitrary line bundle $L$
over a manifold with boundary and extend Demailly's inequalities to
this situation. We will write $h^{q}(L^{k})$ for the dimension of
$H^{0,q}(X,L^{k}),$ the Dolbeault cohomology group of $(0,q)-$forms
with values in $L^{k}.$ The cohomology groups are defined with respect
to forms that are smooth up to the boundary. Recall that $X(q)$ is
the subset of $X$ where $\textrm{index$($}\Theta)=q$ and we let
\[
T(q)_{\rho,x}=\left\{ t>0:\,\textrm{index$($}\Theta+t\mathcal{L})=q\:\textrm{along }T^{1,0}(\partial X)_{x}\right\} .\]
The main theorem we will prove is the following generalization of
Demailly's weak holomorphic Morse inequalities. 

\begin{thm}
\label{thm main} Suppose that $X$ is is a compact manifold with
boundary, such that the Levi form is non-degenerate on the boundary.
Then, up to terms of order $o(k^{n}),$ \begin{equation}
h^{q}(L^{k})\leq k^{n}(-1)^{q}(\frac{1}{2\pi})^{n}(\int_{X(q)}\Theta_{n}+\int_{\partial X}\int_{T(q)_{\rho,x}}(\Theta+t\mathcal{L})_{n-1}\wedge\partial\rho\wedge dt),\label{main statement}\end{equation}

\end{thm}
The integral over $t$ in the boundary integral is oriented so as
to make $(-1)^{q}$ times the last integral non-negative. The boundary
integral above may also be expressed more directly in terms of symplectic
geometry as \begin{equation}
\int_{X_{+}(q)}(\Theta+d\alpha)_{n}\label{intro contact formula}\end{equation}
 where $(X_{+},d\alpha)$ is the so called symplectification of the
contact manifold $\partial X$ induced by the complex structure of
$X$ (section \ref{section hole and contact}).

Examples will be presented that show that the leading constants in
the bounds of the theorem are sharp. We will also obtain the corresponding
generalization of the strong holomorphic Morse inequalities. The most
interesting case is when the manifold is a strongly pseudoconcave
manifold $X$ of dimension $n\geq3$ with a positive line bundle $L.$
Then, if the curvature forms of $L$ and $\partial X$ are conformally
equivalent along the complex tangential directions of $\partial X,$
we will deduce that \begin{equation}
h^{0}(L^{k})=k^{n}(\int_{X}\Theta_{n}+\frac{1}{n}\int_{\partial X}(i\partial\overline{\partial}\rho)_{n-1}\wedge i\partial\rho))+o(k^{n}),\label{intro dim for conf equiv}\end{equation}
 if the defining function $\rho$ is chosen in an appropriate way.
In particular, such a line bundle $L$ is \emph{big} and \ref{intro dim for conf equiv}
can be expressed as \[
\textrm{Vol$(L)=$Vol$(X)+\frac{1}{n}$Vol$(\partial X)$}\]
in terms of the corresponding symplectic volume of $X$ and contact
volume of $\partial X.$ Examples are provided that show that theorem
\ref{thm main} is sharp and also compatible with hole filling

The proof of theorem \ref{thm main} will follow from local estimates
for the corresponding Bergman function $B_{X}^{q,k}$ where $B_{X}^{q,k}$
is the Bergman function of the space $\mathcal{H}^{0,q}(X,L^{k})$
of $\overline{\partial}-$harmonic $(0,q)-$forms satisfying $\overline{\partial}-$Neumann
boundary conditions (simply referred to as the harmonic forms in the
sequel). The point is that the integral of the Bergman function is
the dimension of $\mathcal{H}^{0,q}(X,L^{k})$. It is shown that,
for large $k,$ the Bergman function is estimated by the sum of two
model Bergman functions, giving rise to the bulk and the boundary
integrals in theorem \ref{thm main}. The model at a point $x$ in
the interior of $X$ is obtained by replacing the manifold $X$ with
flat $\C^{n}$ and the line bundle $L$ with the constant curvature
line bundle over $\C^{n}$ obtained by freezing the curvature of the
line bundle at the point $x.$ Similarly, the model at a boundary
point is obtained by replacing $X$ with the unbounded domain $X_{0}$
in $\C^{n},$ whose constant Levi curvature is obtained by freezing
the Levi curvature at the boundary point in $X.$ The line bundle
$L$ is replaced by the constant curvature line bundle over $X_{0},$
obtained by freezing the curvature along the complex tangential directions,
while making it flat in the complex normal direction. 

The method of proof is an elaboration of the, comparatively elementary,
technique introduced in \cite{berm} to handle Demailly's case of
a manifold without boundary. 

\begin{rem}
\label{rem: cond z} The boundary integral in \ref{main statement}
is finite precisely when there is no point in the boundary where the
Levi form $i\partial\overline{\partial}\rho$ has exactly $q$ negative
eigenvalues. Indeed, any sufficiently large $t$ will then be in the
complement of the set $T(q)_{\rho,x}.$ Since, we have assumed that
the Levi form $i\partial\overline{\partial}\rho$ is non-degenerate,
this condition coincides with the so called condition $Z(q)$ \cite{f-k}.
However, for an arbitrary Levi form the latter condition is slightly
more general: it holds if the Levi form has at least $q+1$ negative
eigen values or at least $n-q$ positive eigen values everywhere on
$\partial X.$ In fact, the proof of theorem \ref{thm main} only
uses that $\partial X$ satisfies condition $Z(q)$ and is hence slightly
more general than stated. Furthermore, a function $\rho$ is said
to satisfy condition $Z(q)$ at a point $x$ if $x$ is not a critical
point of $\rho$ and if $\partial\overline{\partial}\rho$ satisfies
the curvature condition at $x$ along the level surface of $\rho$
passing through $x.$ 
\end{rem}
One final remark about the extension of the Morse inequalities to
open manifolds:

\begin{rem}
\label{rem:open cohom}The cohomoloy groups $H^{0,*}(X,L^{k})$ associated
to the manifold with boundary $X$ occuring in the weak Morse inequalities,
theorem \ref{II thm morse}, are defined with respect to forms that
are smooth up to the boundary. Removing the boundary from $X$ we
get an open manifold, that we denote by $\dot{X.}$ By the Dolbeault
theorem \cite{gri} the usual Dolbeault cohomology groups $H^{0,*}(\dot{X},L^{k})$
of $\dot{X}$ are isomorphic to the cohomology groups $H^{*}(\dot{X},\mathfrak{\mathcal{O}}(L^{k}))$
of the sheaf $\mathfrak{\mathcal{O}}(L^{k})$ of germs of holomorphic
sections on $\dot{X}$ with values in $L^{k}.$ Moreover, if we assume
that condition $Z(q)$ \emph{and} $Z(q+1)$ hold then $H^{0,q}(X,L^{k})$
and $H^{0,q}(\dot{X},L^{k})$ are isomorphic \cite{f-k}. Furthermore,
consider a given \emph{open} manifold $Y$ with a smooth exhaustion
function $\rho,$ i.e a function such that the open sublevel sets
of $\textrm{$\rho$}$ are relatively compact in $Y$ for every real
number $c.$ Then, if for a fixed regular value $c_{0},$ the curvature
conditions $Z(q)$ and $Z(q+1)$ hold for $\rho$ when $\rho\geq c_{0},$
the group $H^{0,q}(Y,L^{k})$ is isomorphic to $H^{0,q}(X_{c_{0}},L^{k})$\cite{h},
where $X_{c_{0}}$ is the corresponding closed sublevel set of $\rho.$
In this way one gets Morse inequalities on certain open manifolds
$Y.$ 
\end{rem}
\begin{notation}
The notation $a_{k}\sim(\lesssim)b_{k}$ will stand for $a_{k}=(\leq)C_{k}b_{k},$
where $C_{k}$ tends to one when $k$ tends to infinity. The $\overline{\partial}-$Laplacian
\cite{gri} will be called just the Laplacian. It is the differential
operator defined by $\Delta:=\overline{\partial}\overline{\partial}^{*}+\overline{\partial}^{*}\overline{\partial}$
(where $\overline{\partial}^{*}$ denotes the formal adjoint of $\overline{\partial}$)
acting on smooth forms on $X$ with values in $L^{k}.$ Similarly,
we will call an element in the kernel of $\Delta$ harmonic, instead
of $\overline{\partial}-$harmonic.
\end{notation}
The paper is organized in two parts. In the first part we will state
and prove the \emph{weak} holomorphic Morse inequalites. In the second
part the \emph{strong} holomorphic Morse inequalities are obtained.
Finally, the weak Morse inequalities are shown to be sharp and the
relation to hole filling investigated. 

\tableofcontents{}

\part{\label{II section harmonic}The weak Morse inequalities}

\section{Setup and a sketch of the proof}

In the first part we will show how to obtain weak holomorphic Morse
inequalities for $(0,q)-$ forms, with values in a given line bundle
$L$ over a manifold with boundary $X.$ In other words we will estimate
the dimension of $H^{0,q}(X,L^{k})$ in terms of the curvature of
$L$ and the Levi curvature of the boundary of $X.$ With notation
as in the introduction of the article the theorem we will prove is
as follows.

\begin{thm}
\label{II thm morse}Suppose that $X$ is is a compact manifold with
boundary, such that the Levi form is non-degenerate on the boundary.
Then, up to terms of order $o(k^{n}),$ \[
h^{q}(L^{k})\leq k^{n}(-1)^{q}(\frac{1}{2\pi})^{n}(\int_{X(q)}\Theta_{n}+\int_{\partial X}\int_{T(q)_{\rho,x}}(\Theta+t\mathcal{L})_{n-1}\wedge\partial\rho\wedge dt),\]
(the integral over $t$ is oriented so as to make the last integral
non-negative). 
\end{thm}
Note that the last integral is independent of the choice of defining
function. Indeed, if $\rho'=f\rho$ is another defining function,
where $f$ is a positive function, the change of variables $s=ft$
shows that the integral is unchanged. A more intrinsic formulation
of the last integral will be given in section \ref{section hole and contact}.
Let us now fix the grade $q.$ Since, the statement of the theorem
is vacuous if the Levi form $i\partial\overline{\partial}\rho$ has
exactly $q$ negative eigenvalues somewhere on $\partial X$ (compare
remark \ref{rem: cond z}) we may assume that this is not the case.
Then it is well-known that the Dolbeault cohomology group $H^{0,q}(X,E)$
is finite dimensional for any given vector bundle $E$ over $X.$
The cohomology groups are defined using forms that are smooth up to
the boundary. Moreover, the Hodge theorem, in this context, says that
$H^{0,q}(X,E)$ is isomorphic to the space $\mathcal{H}^{0,q}(X,E),$
consisting of harmonic $(0,q)-$forms, that are smooth up to the boundary,
where they satisfy $\overline{\partial}-$Neumann boundary conditions
(\cite{f-k}). The space $\mathcal{H}^{0,q}(X,E)$ is defined with
respect to a given metrics on $X$ and $E.$

The starting point of the proof of theorem \ref{II thm morse} is
the fact that the dimension of the space $\mathcal{H}^{0,q}(X,L^{k})$
may be expressed as an integral over $X$ of the so called Bergman
function $B_{X}^{q,k}$ defined as\[
B_{X}^{q,k}(x)=\sum_{i}\left|\Psi_{i}(x)\right|^{2},\]
 where $\left\{ \Psi_{i}\right\} $ is any orthonormal base for $\mathcal{H}^{0,q}(X,L^{k}).$
Indeed, the integral of each term in the sum is equal to one. Note
that the Bergman function $B_{X}^{q,k}$ depends on the metric $\omega$
on $X.$ It is convenient to use $k\omega$ as metric for a given
$k.$ The point is that the volume of $X$ mesured with respect to
$k\omega$ is of the order $k^{n}.$ Hence, the dimension bound in
theorem \ref{II thm morse} will follow from a point-wise estimate
of the corresponding Bergman function $B_{X}^{q,k}.$ Another reason
why $k\omega$ is a natural metric on $X$ is that, since $k\phi$
is the induced fiber metric on $L^{k},$ the norms of forms on $X$
with values in $L^{k}$ become more symmetrical with respect to the
base and fiber metrics. In fact, we will have to let the metric $\omega$
itself depend on $k$ (and on a large parameter $R)$ close to the
boundary and we will estimate the Bergman function of the space $\mathcal{H}^{0,q}(X,k\omega_{k},L^{k})$
in terms of model Bergman functions and compute the model cases explicitely.
The sequence of metrics $\omega_{k}$ will be of the following form.
First split the manifold $X$ in an \emph{inner region} $X_{\varepsilon},$
with defining function $\rho+\varepsilon$ and its complement, the
\emph{boundary region,} given a small positive number $\epsilon.$
The level sets where $\rho=-Rk^{-1}$ and $\rho=-k^{-1/2}$ divide
the boundary region into three regions. The one that is closest to
the boundary of $X$ will be called the \emph{first region} and so
on. Next, define $\omega_{T},$ the complex tangential part of $\omega$
close to the boundary by \[
\omega_{T}:=\omega-2i\partial\rho\wedge\overline{\partial\rho}\]
 (recall that we assumed that $d\rho$ is of unit-norm with respect
to $\omega$ close to the boundary of $X).$ The metric $\omega_{k}$
is of the form \begin{equation}
\omega_{k}=\omega_{T}+a_{k}(\rho)^{-1}2i\partial\rho\wedge\overline{\partial\rho},\label{II metric express}\end{equation}
 where the sequence of smooth functions $a_{k}$ will be chosen so
that, basically, the distance to the boundary, when mesured with respect
to $k\omega_{k},$ in the three different regions is independent of
$k.$ The properties of $\omega_{k}$ that we will use in the two
regions will be stated in the proofs below, while the precise definition
of $\omega_{k}$ is deferred to section \ref{II subsection metricsequence}.

\subsection{\label{sec:A-sketch-of}A sketch of the proof of the weak Morse inequalities }

To make the sketch of the proof cleaner, we will just show how to
estimate the \emph{extremal function} \begin{equation}
S_{X}^{q,k}(x)=\sup_{\alpha_{k}}\left|\alpha_{k}(x)\right|^{2}\label{II sketch def of s}\end{equation}
 closely related to $B_{X}^{q,k},$ where the supremum is taken over
all normalized elements of the space $\mathcal{H}^{0,q}(X,\omega_{k},L^{k}).$
When $q=0,$ i.e the case of holomorphic sections, it is a classical
fact that they are actually equal and the general relation is given
in section \ref{II section bergman}. Let us first see how to get
the following bound in the inner region $X_{\varepsilon}$ defined
above: \begin{equation}
S_{X}^{q,k}(x)\lesssim S_{\C^{n},x}^{q}(0),\label{II sketch bound on s inner}\end{equation}
where the right hand side is the extremal function for the model case
defined below. Moreover, the left hand side is uniformly bounded by
a constant, which is essential when integrating the estimate to get
an estimate on the dimension of $\mathcal{H}^{0,q}(X,k\omega_{k},L^{k}).$
The proof of \ref{II sketch bound on s inner} proceeds exactly as
in the case when $X$ is a compact manifold without boundary \cite{berm}.
Let us recall the argument, slightly reformulated. Fix a point $x$
in $X_{\varepsilon}.$ We may take local holomorphic coordinates centered
at $x$ and a local trivialization of $L$ such that\begin{equation}
\phi(z)=\sum_{i=1}^{n}\lambda_{i}z_{i}\overline{z_{i}}+...,\,\,\,\omega(z)=\frac{i}{2}\sum_{i=1}^{n}dz_{i}\wedge\overline{dz_{i}}+....\label{eq: sktech metrics}\end{equation}
 where the dots indicate lower order terms and the leading terms are
called model metrics and denoted by $\phi_{0}$ and $\omega_{0},$
respectively. \emph{}Hence, the model situation is a line bundle of
constant curvature on flat $\C^{n}.$ Note that the unit ball at $x$
with respect to the metric $k\omega_{k}$ corresponds approximately
to the coordinate ball at $0$ of radius $k^{-1/2}.$ To make this
more precise, define a scaling map \[
F_{k}(z)=k^{-1/2}z\]
 and consider a sequence of expanding balls centered at $0$ in $\C^{n}$
of radius $r_{k},$ slowly exhausting all of $\C^{n}.$ We will call
$F_{k}^{*}(k\phi)$ and $F_{k}^{*}(k\omega)$ the scaled metrics on
the expanding balls. The point is that they converge to the model
metrics $\phi_{0}$ and $\omega_{0}.$ This follows immediately from
the expression \ref{eq: sktech metrics} and the fact that the model
metrics are invariant when $kF_{k}^{*}$ is applied. Next, given a
$(0,q)$ form on $X$ with values in $L^{k},$ we denote by $\alpha^{(k)}$
the scaled form defined by $\alpha^{(k)}=F_{k}^{*}\alpha_{k}.$ Then,
by the convergence of the scaled metrics,\begin{equation}
\left\Vert \alpha_{k}\right\Vert _{F_{k}(B_{r_{k}})}^{2}\sim\left\Vert \alpha^{(k)}\right\Vert _{B_{r_{k}}}^{2},\label{eq: local of norms}\end{equation}
using the norms induced by the model metrics in the right hand side
above. Now, if $\alpha_{k}$ is a normalized sequence of extremals
(i.e realizing the extremum in \ref{II sketch def of s}) we have\[
S_{X}^{q,k}(x)=\left|\alpha^{(k)}(0)\right|^{2}.\]
 By \ref{eq: local of norms}, the norms of the scaled sequence $\alpha^{(k)}$
are less than one, when $k$ tends to infinity. Moreover, $\alpha^{(k)}$
is harmonic with respect to the scaled metrics and since these converge
to model metrics, inner elliptic estimates for the Laplacian show
that there is a subsequence of $\alpha^{(k)}$ that converges to a
model harmonic form $\beta$ in $\C^{n}.$ In fact, we may assume
that the whole sequence $\alpha^{(k)}$ converges. Hence, \[
\limsup_{k}\left|\alpha^{(k)}(0)\right|^{2}=\left|\beta(0)\right|^{2}\]
 which in turn is bounded by the model extremal function $S_{X_{0},x}(0).$
Moreover, since $X_{\varepsilon}$ may be covered by coordinate balls
of radius $k^{-1/2},$ staying inside of $X$ for large $k,$ one
actually gets a uniform bound. 

Let us now move on to the boundary region $X-X_{\varepsilon}$ that
we split into three regions as above. Fix a point $\sigma$ in the
boundary of $X.$ We may take local coordinates centered at $\sigma$
and orthonormal at $\sigma,$ so that \[
\rho(z,w)=v-\sum_{i=1}^{n-1}\mu_{i}\left|z_{i}\right|^{2}+...\]
 where $v$ is the imaginary part of $w$ \cite{cm}. The leading
term of $\rho$ will be denoted by $\rho_{0},$ and will be referred
to as the defining function of the \emph{model domain} $X_{0}$ in
$\C^{n}.$ Observe that the model domain $X_{0}$ is invariant under
the holomorphic anisotropic scaling map \[
F_{k}(z,w)=(z/k^{1/2},w/k).\]
Moreover, the scaled fiber metric on $L^{k}$ now tends to the new
model fiber metric $\phi_{0}(z,0),$ since the terms in $\phi$ involving
the coordinate $w$ are suppressed by the anisotropic scaling map
$F_{k}.$ Now, the bound \ref{II sketch bound on s inner} is replaced
by \begin{equation}
S_{X}^{q,k}(0,iv/k)\lesssim S(0,iv),\label{II sketch s bound in boundary}\end{equation}
 in terms of the new model case. To see this one replaces the balls
of decreasing radii used before with $F_{k}(D_{k})$ intersected with
$X,$ where $D_{k}$ is a sequence of slowly expanding polydiscs.
Moreover, we have to let the initial metric $\omega$ on $X$ depend
on $k$ in the normal direction in order that the scaled metric converge
to a non-degenerate model metric. In the first region we will essentially
let\[
\omega_{k}=\omega_{T}+k2i\partial\rho\wedge\overline{\partial\rho}.\]
 As a \emph{model metric} in $X_{0}$ we will essentially use \[
\omega_{0}=\frac{i}{2}\partial\overline{\partial}\left|z\right|^{2}+2i\partial\rho_{0}\wedge\overline{\partial\rho_{0}}.\]
 Then clearly \begin{equation}
F_{k}^{*}(k\omega_{k})=\omega_{0}\label{II scaling of base metrics}\end{equation}
 in the model case and it also holds asymptotically in $k,$ in the
general case. Replacing the inner elliptic estimates used in the inner
part $X_{\varepsilon}$ with subelliptic estimates for the $\overline{\partial}-$Laplacian
close to the boundary one gets the bound \ref{II sketch s bound in boundary}
more or less as before. Finally, using similar scaling arguments,
one shows that the contribution from the second and third region to
the total integral of $B_{X}^{q,k}$ is negligible when $k$ tends
to infinity. 

This gives the bound \[
\int_{X}B_{X}^{q,k}(k\omega_{k})_{n}\lesssim k^{n}(\int_{X}B_{\C^{n},x}^{q}\omega_{n}+\int_{\partial X}\int_{-\infty}^{0}B_{X_{0},\sigma}^{q}(iv)dvd\sigma)\]
 integrating over an infinite ray in the model region $X_{0}$ in
the second integral (after letting $R$ tend to infinity). Computing
the model Bergman functions explicitly then finishes the proof of
the theorem.

\section{\label{II section bergman}Bergman kernel forms}

Let us now turn to the detailed proof of theorem \ref{II thm morse}.
First we introduce Bergman kernel forms to relate the Bergman function
$B_{X}^{q,k}$ to extremal functions taking account of the components
of a form (see \cite{berm2} for proofs). Let $(\psi_{i})$ be an
orthonormal base for a finite dimensional Hilbert space $\mathcal{H}^{0,q}$
of $(0,q)-$ forms with values in $L.$ Denote by $\pi_{1}$ and $\pi_{2}$
the projections on the factors of $X\times X.$ The \emph{Bergman
kernel} \emph{form} of the Hilbert space $\mathcal{H}^{0,q}$ is defined
by \[
\K^{q}(x,y)=\sum_{i}\psi_{i}(x)\wedge\overline{\psi_{i}(y)}.\]
 Hence, $\K^{q}(x,y)$ is a form on $X\times X$ with values in the
pulled back line bundle $\pi_{1}^{*}(L)\otimes\overline{\pi_{2}^{*}(L)}.$
For a fixed point $y$ we identify $\K_{y}^{q}(x):=\K^{q}(x,y)$ with
a $(0,q)-$form with values in $L\otimes\Lambda^{0,q}(X,\overline{L})_{y}.$
The definition of $\K^{q}$ is made so that $\K^{q}$ satisfies the
following reproducing property: \begin{equation}
\alpha(y)=c_{n,q}\int_{X}\alpha\wedge\overline{\K_{y}^{q}}\wedge e^{-\phi}\omega_{n-q},\label{II repr prop}\end{equation}
for any element $\alpha$ in $\mathcal{H}^{0,q},$ using a suggestive
notation and where $c_{n,q}$ is a complex number of unit norm that
ensures that \ref{II repr prop} may be interpreted as a scalar product.
Properly speaking, $\alpha(y)$ is equal to the push forward $\pi_{2*}(c_{n,q}\alpha\wedge\K^{q}\wedge\omega_{n-q}e^{-\phi})(y).$
The restriction of $\K^{q}$ to the diagonal can be identified with
a $(q,q)-$form on $X$ with values in $L\otimes\overline{L}.$ The
\emph{Bergman form} is defined as $\K^{q}(x,x)e^{-\phi(x)},$ i.e.
\begin{equation}
\B^{q}(x)=\sum_{i}\psi_{i}(x)\wedge\overline{\psi_{i}(x)}e^{-\phi(x)}\label{II eq def of b form}\end{equation}
and it is a globally well-defined $(q,q)-$form on $X.$ The following
notation will turn out to be useful. For a given form $\alpha$ in
$\Omega^{0,q}(X,L)$ and a decomposable form in $\Omega^{0,q}(X)_{x}$
of unit norm, let $\alpha_{\theta}(x)$ denote the element of $\Omega^{0,0}(X,L)_{x}$
defined as \[
\alpha_{\theta}(x)=\left\langle \alpha,\theta\right\rangle _{x}\]
where the scalar product takes values in $L_{x}.$ We call $\alpha_{\theta}(x)$
the \emph{value of $\alpha$ at the point $x,$ in the direction $\theta.$}
Similarly, let $B_{\theta}^{q}(x)$ denote the function obtained by
replacing \ref{II eq def of b form} by the sum of the squared pointwise
norms of $\psi_{i,\theta}(x)$. Then $B_{\theta}^{q}(x)$ has the
following useful extremal property:\begin{equation}
B_{\theta}^{q}(x)=\sup_{\alpha}\left|\alpha_{\theta}(x)\right|^{2},\label{II extremal prop of b direction}\end{equation}
where the supremum is taken over all elements $\alpha$ in $\mathcal{H}^{0,q}$
of unit norm. The supremum will be denoted by $S_{\theta}^{q}(x)$
and an element $\alpha$ realizing the supremum will be referred to
as an \emph{extremal form for the space $\mathcal{H}^{0,q}$ at the
point $x,$ in the direction $\theta.$} The reproducing formula \ref{II repr prop}
may now be written as \[
\alpha_{\theta}(y)=(\alpha,\K_{y,\theta}^{q}).\]
Finally, note that the Bergman function $B$ is the trace of $\B^{q},$
i.e. \[
B^{q}\omega_{n}=c_{n,q}\B^{q}\wedge\omega_{n-q}.\]
Using the extremal characterization \ref{II extremal prop of b direction}
we have the following useful expression for $B:$ \begin{equation}
B(x)=\sum_{\theta}S_{\theta}(x),\label{II extremal prop of b}\end{equation}
 where the sum is taken over any orthonormal base of direction forms
$\theta$ in $\Lambda^{0,q}(X)_{x}.$

\section{\label{II section model}The model boundary case}

In the sketch of the proof of the weak holomorphic Morse inequalities
(section \ref{sec:A-sketch-of}) it was explained how to bound the
Bergman function on $X$ by model Bergman functions. In this section
we will compute the Bergman kernel explicitely in the model boundary
case. Consider $\C^{n}$ with coordinates $(z,w),$ where $z$ is
in $\C^{n-1}$ and $w=u+iv.$ Let $X_{0}$ be the domain with defining
function \[
\rho_{0}(\textrm{z},w)=v+\psi_{0}(z):=v+\sum_{i=1}^{n-1}\mu_{i}\left|z_{i}\right|^{2},\]
 and with the metric \[
\omega_{0}=\frac{i}{2}\partial\overline{\partial}\left|z\right|^{2}+a(\rho)^{-1}2i\partial\rho_{0}\wedge\overline{\partial\rho_{0}}.\]
Note that the corresponding volume element $(\omega_{0})_{n}$ is
given by $a(\rho)^{-1}$ times the usual Euclidian volume element
on $\C^{n}.$ We will take $a(\rho_{0})$ to be comparable to $(1-\rho_{0})^{2}$
(compare section \ref{II subsection metricsequence}) but we will
only use that the corresponding metric $\omega_{0}$ is {}``relatively
complete'' (compare section \ref{II subsection proof of tangent}).
We fix the $q$ and assume that condition $Z(q)$ holds on $\partial X_{0},$
i.e. that at least $q+1$ of the eigen values $\mu_{i}$ are negative
or that at least $n-q$ of them are positive. 

Let $\mathcal{H}^{0,q}(X_{0},\phi_{0})$ be the space of all $(0,q)-$
forms on $X$ that have finite $L^{2}-$ norm with respect to the
norms defined by the metric $\omega_{0}$ and the weight $e^{-\phi_{0}(z)},$
where $\phi_{0}$ is quadratic, and that are harmonic with respect
to the corresponding Laplacian.%
\footnote{Since the metric is {}``relatively complete'', all harmonic forms
are in fact closed and coclosed with respect to $\overline{\partial}.$
One could also take this as the definition of harmonic forms in the
present context. %
} Moreover, we assume the the forms are smooth up to the boundary of
$X_{0}$ where they satisfy $\overline{\partial}-$Neumann boundary
conditions (the regularity properties are automatic, since we have
assumed that condition $Z(q)$ holds \cite{f-k}). The Bergman kernel
form of the Hilbert space $\mathcal{H}^{0,q}(X_{0},\phi_{0})$ will
be denoted by $\K_{X_{0}}^{q}.$ We will show how to expand $\K_{X_{0}}^{q}$
in terms of Bergman kernels on $\C^{n-1},$ and then compute it explicitly.
Note that the metric $\omega_{0}$ is chosen so that the pullback
of any form on $\C^{n-1}$ satisfies $\overline{\partial}-$Neumann
boundary conditions. Conversely, we will show that any form in $\mathcal{H}^{0,q}(X_{0},\phi_{0})$
can be written as a superposition of such pulled-back forms. 

By the very definition of the metric $\omega_{0},$ the forms $dz_{i}$
and $a^{-1/2}\partial\rho_{0}$ together define an orthogonal frame
of $(1,0)-$ forms. Any $(0,q)-$ form $\alpha$ on $X$ may now be
uniquely decomposed in a {}``tangential'' and a {}``normal'' part:
\[
\alpha=\alpha_{T}+\alpha_{N},\]
 where $\alpha=\alpha_{T}$ modulo the algebra generated by $\overline{\partial\rho_{0}}.$
A form $\alpha$ without normal part will be called \emph{tangential.}
The proof of the following proposition is postponed till the end of
the section.

\begin{prop}
\label{II prop tangential}Suppose that $\alpha$ is in $\mathcal{H}^{0,q}(X_{0},\phi_{0}).$
Then $\alpha$ is tangential, closed and coclosed (with respect to
$\overline{\partial}).$
\end{prop}
By the previous proposition any form $\alpha$ in $\mathcal{H}^{0,q}(X_{0},\phi_{0})$
may be written as \[
\alpha(z,w)=\sum_{I}f_{I}d\overline{z_{I}}\]
Moreover, since $\alpha$ is in $L^{2}(X_{0},\phi_{0})$ and $\overline{\partial}-$
closed, the components $f_{I}$ are in $L^{2}(X_{0},\phi_{0})$ and
holomorphic in the $w-$ variable. We will have use for the following
basic lemma:%
\footnote{This lemma can be reduced to the Payley-Wiener theorem 19.2 in \cite{ru}. %
}

\begin{lem}
Let $m(v)$ be a positive function on $[0,\infty[$ with polynomial
growth at infinity. If $f(w)$ is a a holomorphic function in $\left\{ v<c\right\} $
with finite $L^{2}-$norm with respect to the measure $m(v)dudv,$
then there exists a function $\widehat{f(t)}$ on $]0,\infty[$ such
that \[
f(w)=\int_{0}^{\infty}\widehat{f(t)}e^{-\frac{i}{2}wt}dt\]
 Moreover, \begin{equation}
\int_{v<c}\int_{u=-\infty}^{\infty}\left|f(w)\right|^{2}m(v)dudv=4\pi\int_{v<c}\int_{t=0}^{\infty}\left|\widehat{f}(t)\right|^{2}e^{vt}m(v)dtdv\label{eq: parcev}\end{equation}

\end{lem}
We will call $\widehat{f(t)}$ the Fourier transform of $f(w).$ Now,
fix $z$ in $\C^{n-1}$ and take $c=-\psi_{0}(z)$ and $m(v)=a(\rho_{0})^{-1}=a(v+\psi_{0}(z))^{-1}.$
Then $f_{I},$ as a function of $w,$ must satisfy the requirements
in the lemma above for almost all $z.$ Fixing such a $z$ we write
$\widehat{f_{I,t}(z)}$ for the function of $t$ obtained by taking
the Fourier transformation with respect to $w.$

Hence, we can write \begin{equation}
\alpha(z,w)=\int_{0}^{\infty}\widehat{\alpha_{t}}(z)e^{-\frac{i}{2}wt}dt\label{eq: fouriertransf for alfa}\end{equation}
for almost all $z,$ if we extend the Fourier transform and the integral
to act on forms coefficientwise. Note that the equality \ref{eq: fouriertransf for alfa}
holds in $L^{2}(X_{0},\phi_{0}).$ The following proposition describes
the space $\mathcal{H}^{0,q}(X_{0},\phi_{0})$ in terms of the spaces
$\mathcal{H}^{0,q}(\C^{n-1},t\psi_{0}+\phi_{0}),$ consisting of all
harmonic $(0,q)-$ forms in $L^{2}(\C^{n-1},t\psi_{0}+\phi_{0})$
(with respect to the Euclidean metric in $\C^{n}).$ The corresponding
scalar products over $\C^{n-1}$ are denoted by $(\cdot,\cdot)_{t}.$ 

\begin{prop}
\label{II lemma model scalar pr}Suppose that $\alpha$ is a tangential
$(0,q)-$form on $X_{0}$ with coefficients holomorphic with respect
to $w.$ Then $\widehat{\alpha_{t}}$ is in $L^{2}(\C^{n-1},t\psi_{0}+\phi_{0})$
for almost all $t$ and \begin{equation}
(\alpha,\alpha)_{X_{0}}=4\pi\int(\widehat{\alpha_{t}},\widehat{\alpha_{t}})_{t}b(t)dt,\label{II prop scalar pr statement}\end{equation}
 where $b(t)=\int_{s<0}e^{st}a(\rho_{0})^{-1}ds.$ Moreover, if $\alpha$
is in $\mathcal{H}^{0,q}(X_{0},\phi_{0}),$ then $\widehat{\alpha_{t}}$
is in $\mathcal{H}^{0,q}(\C^{n-1},t\psi_{0}+\phi_{0})$ for almost
all $t.$
\end{prop}
\begin{proof}
It is clearly enough to prove \ref{II prop scalar pr statement} for
the components $f_{I}$ of $\alpha,$ i.e for a function $f$ in $X_{0}$
that is holomorphic with respect to $w.$ When evaluating the norm
$(f,f)_{X_{0}}$ over $X_{0}$ we may first perform the integration
over $u,$ using \ref{eq: parcev}, giving\[
(f,f)_{X_{0}}=4\pi\int_{\rho_{0}<0}\left|\widehat{f_{t}}(z)\right|^{2}e^{vt}e^{-\phi(z)}a(\rho_{0})^{-1}dzdtdv,\]
 where $dz$ stands for the Euclidean volume form $(\frac{i}{2}\partial\overline{\partial}\left|z\right|^{2})_{n-1}$
on $\C^{n-1}.$ If we now fix $z$ and make the change of variables
$s:=v+\psi_{0}(z)$ and integrate with respect to $s$ we get\[
4\pi\int\left|\widehat{f_{t}}(z)\right|^{2}b(t)e^{-(t\psi(z)+\phi(z))}dtdz=4\pi\int(\widehat{f_{t}},\widehat{f}_{t})_{t}b(t)dt.\]
 Since this integral is finite, it follows that $(\widehat{f_{t}},\widehat{f_{t}})$
is finite for almost all $t.$

Next, assume that $\alpha$ is in $\mathcal{H}^{0,q}(X_{0},\phi_{0}).$
By proposition \ref{II prop tangential} $\alpha$ is $\overline{\partial}-$closed,
so that \ref{eq: fouriertransf for alfa} gives that $\widehat{\alpha_{t}}$
is $\overline{\partial}-$closed for almost all $t.$ Let us now show
that $\widehat{\alpha_{t}}$ is $\overline{\partial}-$coclosed with
respect to $L^{2}(\C^{n-1},t\psi_{0}+\phi_{0})$ for almost all $t.$
Fix an interval $I$ in the positive half-line and let $\beta$ be
a form in $X_{0}$ that can be written as \[
\beta(z,w)=\int_{t\ \in I}\eta_{t}(z)e^{-\frac{i}{2}wt}dt\]
where $\eta_{t}$ is a smooth $(0,q-1)-$form with compact support
on $\C^{n-1}$ for a fixed $t$ (and mesurable with respect to $t$
for $z$ fixed). In particular $\beta$ is a smooth form in $L^{2}(X_{0},\phi_{0})$
that is tangential and holomorphic with respect to $w.$ According
to \ref{eq: fouriertransf for alfa} $\widehat{\beta_{t}}$ is equal
to $\eta_{t}$ for $t\in I$ and vanishes otherwise. By proposition
\ref{II prop tangential} $\alpha$ is $\overline{\partial}-$coclosed
(with respect to $L^{2}(X_{0},\phi_{0})).$ Using \ref{II prop scalar pr statement}
we get that \[
0=(\overline{\partial}^{*}\alpha,\beta)=(\alpha,\overline{\partial}\beta)=4\pi\int_{t\in I}(\widehat{\alpha_{t}},\overline{\partial}\eta_{t})_{t}b(t)dt,\]
where we have used lemma \ref{II formal adj} proved in the next section
to get the second equality. Since this holds for any choice of form
$\beta$ and interval $I$ as above we conclude that $\overline{\partial}^{*}\widehat{\alpha_{t}}=0$
for almost all $t.$ Hence $\widehat{\alpha_{t}}$ is in $\mathcal{H}^{0,q}(\C^{n-1},t\Psi_{0}+\phi_{0})$
for almost all $t.$
\end{proof}
Denote by $\K_{t}^{q}$ the Bergman kernel of the Hilbert space $\mathcal{H}^{0,q}(\C^{n-1},t\Psi_{0}+\phi_{0}).$ 

\begin{lem}
\label{II lemma model k expansion} The Bergman kernel $\K_{X_{0}}^{q}$
may be expressed as \[
\K_{X_{0}}^{q}(z,w,z',w')=\frac{1}{4\pi}\int_{0}^{\infty}\K_{t}^{q}(z,z')e^{\frac{i}{2}(\overline{w}-w')t}b(t)^{-1}dt\]
In particular, the Bergman form $\B_{X_{0}}^{q}$ is given by \[
\B_{X_{0}}^{q}(z,w)=\frac{1}{4\pi}\int_{0}^{\infty}\B_{t}^{q}(z,z)e^{\rho_{0}t}b(t)^{-1}dt.\]

\end{lem}
\begin{proof}
Take a form $\alpha$ is in $\mathcal{H}^{0,q}(X_{0},\phi_{0})$ and
expand it in terms of its Fourier transform as in \ref{eq: fouriertransf for alfa}.
According to the previous lemma $\widehat{\alpha{}_{t}}$ is in $\mathcal{H}^{0,q}(\C^{n-1},t\Psi_{0}+\phi_{0})$
for almost all $t.$ Hence, we can express it in terms of the corresponding
Bergman kernel $\K_{t}^{q},$ giving \[
f_{I}(z,w)=\int\widehat{f_{I,t}}(z)e^{-\frac{i}{2}wt}dt=\int(\widehat{\alpha_{t}}(z),\K_{t,z,I}^{q})e^{\frac{i}{2}wt}dt,\]
 where $\K_{t,z,I}^{q}$ denotes the Bergman kernel form $\K_{t}^{q}$
at the point $z$ in the direction $d\overline{z_{I}}$ (see section
\ref{II section bergman}) and where have used the reproducing property
\ref{II repr prop} of the Bergman kernel. Now, using the relation
between the different scalar products in the previous lemma we get
\[
f_{I}(z,w)=\frac{1}{4\pi}(\alpha(z',w'),\int_{t}\K_{t,z,I}^{q}(z')e^{\frac{i}{2}\overline{w}t}e^{-\frac{i}{2}w't}b(t)^{-1}dt)_{X_{0}}\]
where $(z',w')$ are the integration variables in the scalar product.
But this means exactly that $\K_{X_{0}}^{q}$ as defined in the statement
of the lemma is the Bergman kernel form of the space $\mathcal{H}^{0,q}(X_{0},\phi_{0})$
since $\alpha$ was chosen arbitrarily. Finally, by definition we
have that $\B_{X_{0}}^{q}(z,w)=\K_{X_{0}}^{q}(z,w,z,w)e^{-\phi_{0}(z)}$
and $\B_{t}^{q}(z,z)=\K_{t}^{q}(z,z')e^{-(t\Psi_{0}+\phi_{0})(z)}.$
Hence, the expression for $\B_{X_{0}}^{q}$ is obtained. 
\end{proof}
Now we can give an explicit expression for the Bergman kernel form
and the Bergman function. In the formulation of the following theorem
we consider $X_{0}$ as a fiber bundle of infinity rays $]-\infty,0]$
over (or rather under) the boundary $\partial X_{0}.$ Then we can
consider the fiber integral over $0$, i.e. the push forward at $0$,
of forms on $X_{0}.$ Moreover, given a real-valued function $\eta$
on $\C^{n-1}$ such that $\frac{i}{2}\partial\overline{\partial}\eta$
has exactly $q$ negative eigenvalues, we define an associated $(q,q)-$
form $\chi^{q,q}$ by \[
\chi^{q,q}:=(i/2)^{q}e^{1}\wedge\cdot\cdot\cdot\wedge e^{q}\wedge\overline{e^{1}}\wedge\cdot\cdot\cdot\wedge\overline{e^{q}}\]
 where $e^{i}$ is a orthonormal $(1,0)-$frame that is dual to a
base $e_{i}$ of the direct sum of eigen spaces corresponding to negative
eigenvalues of $i\partial\overline{\partial}\eta$ (compare \cite{berm2}).
The $(q,q)-$ form associated to $\frac{i}{2}\partial\overline{\partial}\phi_{0}+t\frac{i}{2}\partial\overline{\partial}\rho_{0}$
is denoted by $\chi_{t}^{q,q}$ in the statement of the following
theorem. 

\begin{thm}
\label{II thm model}The Bergman form $\B_{X_{0}}^{q}$ can be written
as an integral over a parameter $t:$ \[
\B_{X_{0}}^{q}(0,u+iv)=\frac{1}{4\pi}\frac{1}{\pi^{n-1}}\int_{T(q)}\chi_{t}^{q,q}\det(\frac{i}{2}\partial\overline{\partial}\phi_{0}+t\frac{i}{2}\partial\overline{\partial}\rho_{0})e^{vt}b(t)^{-1}dt,\]
 where $b(t)=\int_{\rho<0}e^{\rho t}a(\rho)^{-1}d\rho.$ In particular,
the fiber integral over $0$ of the Bergman function $B_{X_{0}}$
times the volume form is given by \begin{equation}
\int_{v=-\infty}^{0}B_{X_{0}}^{q}(0,iv)(\omega_{0})_{n}=(\frac{i}{2\pi})^{n}(-1)^{q}\int_{T(q)}(\partial\overline{\partial}\phi_{0}+t\partial\overline{\partial}\rho_{0})_{n-1}\wedge\partial\rho\wedge dt.\label{II statement of thm model }\end{equation}

\end{thm}
\begin{proof}
Let us first show how to get the expression for $\B_{X_{0}}^{q}.$
Using the previous proposition we just have to observe that in $\C^{m},$
with $\eta$ a quadratic weight function, the Bergman form is given
by \begin{equation}
\B_{\eta}^{q}=\frac{1}{\pi^{m}}1_{X(q)}\chi^{q,q}\det(\frac{i}{2}\partial\overline{\partial}\eta),\label{II pf of theorem model}\end{equation}
 where the constant function $1_{X(q)}$ is equal to one if $\frac{i}{2}\partial\overline{\partial}\eta$
has precisely $q$ negative eigenvalues and is zero otherwise (see
\cite{berm}, \cite{berm2}, \cite{bern}). Next, from section \ref{II section bergman}
we have that $B_{X_{0}}^{q}(\omega_{0})_{n}$ is given by $\B_{X_{0}}^{q}(\omega_{0})_{n-q}.$
Note that \[
\chi^{q,q}\wedge(\omega_{0})_{n-q}=(\frac{i}{2}\partial\overline{\partial}\left|z\right|^{2})_{n-1}\wedge a(\rho_{0})^{-1}(2i)\partial\rho_{0}\wedge d\rho_{0},\]
 Thus, the fiber integral of $B_{X_{0}}^{q}(\omega_{0})_{n}$ over
$0$ reduces to \ref{II statement of thm model } since the factor
$b(t)$ is cancelled by the integral of $e^{vt}a(v)^{-1}.$ 
\end{proof}
Let us finally prove proposition \ref{II prop tangential} .

\subsection{\label{II subsection proof of tangent}The proof of proposition \ref{II prop tangential}:
all harmonic forms are tangential, closed and coclosed.}

We may write \[
\omega_{0}=\frac{i}{2}\partial\overline{\partial}\left|z\right|^{2}+2i\partial\rho'\wedge\overline{\partial\rho'}\]
 for a certain function $\rho'$of $\rho_{0}.$ The forms $2^{-1/2}dz_{i}$
and $2^{1/2}\partial\rho'$ together define a orthonormal frame of
$(1,0)-$ forms. However, we will use the orthogonal frame consisting
of all $dz_{i}$ and $\partial\rho'$ in order not to clutter the
formulas. A dual frame of $(1,0)-$ vector fields is obtained as \begin{equation}
Z_{i}:=\frac{\partial}{\partial z_{i}}-2i\mu_{i}\overline{z_{i}}\frac{\partial}{\partial w},i=1,2,...,n\,\,\,\,\, N:=ia^{1/2}\frac{\partial}{\partial w},\label{II proof of tang vector}\end{equation}
 where $Z_{i}$ is tangential to the level surfaces of $\rho,$ while
$N$ is a complex normal vector field. We decompose any form $\alpha$
as \[
\alpha=\alpha_{T}+\alpha_{N}=\sum f_{I}d\overline{z_{I}}+\overline{\partial\rho}'\wedge g^{0,q-1}\]
Similarly, we decompose the $\overline{\partial}-$ operator acting
on the algebra of forms $\Omega^{0,*}(X_{0})$ as \begin{equation}
\overline{\partial}=\overline{\partial_{T}}+\overline{\partial_{N}}=\sum_{i=1}^{n-1}\overline{Z_{i}}d\overline{z_{i}}\wedge+\overline{N}\overline{\partial\rho'}\wedge,\label{II eq exp of dbar}\end{equation}
 where the vector fields $\overline{Z_{i}}$ etc act on forms over
$X_{0}$ by acting on the coefficients where $d\overline{z_{i}}\wedge$
etc denotes the operator acting on forms on $X,$ obtained by wedging
with $d\overline{z_{i}}.$ The adjoint operator will be denoted by
$d\overline{z_{i}}^{*}.$ Note that the expression for $\overline{\partial}$
is independent of the ordering of the operators, since the elements
in the corresponding frame of $(0,1)-$forms are $\overline{\partial}-$closed.
We denote by $\Delta_{T}$ and $\Delta_{N}$ the corresponding Laplace
operators, i.e. \[
\Delta_{T}=\overline{\partial_{T}}\overline{\partial_{T}}^{*}+\overline{\partial_{T}}^{*}\overline{\partial_{T}},\,\,\,\Delta_{N}=\overline{\partial_{N}}\overline{\partial_{N}}^{*}+\overline{\partial_{N}}^{*}\overline{\partial_{N}},\]
Recall that $\alpha$ is said to satisfy $\overline{\partial}-$Neumann
boundary conditions if $\overline{\partial\rho'}^{*}$applied to $\alpha$
and $\overline{\partial}\alpha$ vanishes on the boundary of $X_{0},$
or equivalently if $\alpha_{N}$ and $(\overline{\partial}\alpha)_{N}$
vanishes there.

\begin{lem}
\label{II formal adj}Denote by $\overline{Z}_{i}^{*}$and $\overline{N}^{*}$the
formal adjoint operators of the operators $\overline{Z}_{i}$and $\overline{N}$
acting on $\Omega^{0,*}(X_{0})$. Then \begin{equation}
\begin{array}{rcl}
\overline{Z_{i}}^{*} & = & -e^{\phi_{0}}Z_{i}e^{-\phi_{0}}\\
\overline{N}^{*} & = & -a^{1/2}Na^{-1/2}\end{array}\label{II lemma formal adj statement}\end{equation}
 Moreover, if the form $\alpha$ has relatively compact support in
$X_{0}$ then for any form smooth form $\beta$ in $X_{0}$ we have
that $\overline{(\partial_{T}}^{*}\alpha,\beta)=(\alpha,\overline{\partial_{T}}\beta)$
and if furthermore $\alpha$ satisfies $\overline{\partial}-$Neumann
boundary conditions, then $\overline{(\partial_{N}}^{*}\alpha,\beta)=(\alpha,\overline{\partial_{N}}\beta)$
(in terms of the formal adjoint operators). 
\end{lem}
\begin{proof}
It is clearly enough to prove \ref{II lemma formal adj statement}
for the action of the operators on smooth \emph{functions} with compact
support (i.e we write $\alpha=f$ and $\beta=g,$ where $f$ and $g$
are smooth functions with compact support). To prove the first statement
in \ref{II lemma formal adj statement} it is, using Leibniz rule,
enough to show that \[
\int_{X}\overline{(Z_{i}}(f\overline{g}e^{-\phi_{0}})(\omega_{0})_{n}=0\]
But this follows from Stokes theorem since the integrand can be written
as a constant times the form\[
d(f\overline{g}e^{-\phi_{0}}a^{-1}(\bigwedge_{j\neq i}dz_{j}\wedge d\overline{z_{j}})\wedge dz_{i}\wedge dw\wedge d\overline{w}),\]
using that $\overline{Z_{i}}(a^{-1})=0,$ since $\overline{Z_{i}}$
is tangential. Similarly, to prove the second statement in \ref{II lemma formal adj statement}
it is enough to observe that \[
\int_{X}d(f\overline{g}e^{-\phi_{0}}a^{-1/2}(\partial\overline{\partial}\left|z\right|^{2})^{n-1}\wedge dw)=0,\]
by Stokes theorem. Indeed, we have that $\overline{N}:=-ia^{1/2}\frac{\partial}{\partial\overline{w}},$
so the statement now follows from Leibniz' rule. Finally, the last
two statements follow from the arguments above, since the boundary
integrals obtained from Stokes theorem vanish.
\end{proof}
\begin{lem}
\label{II lem decomp of lapl}The $\overline{\partial}-$Laplacian
$\Delta$ acting on $\Omega^{0,*}(X_{0})$ decomposes as\[
\Delta=\Delta_{T}+\Delta_{N}.\]

\end{lem}
\begin{proof}
Expanding with respect to the decomposition \ref{II eq exp of dbar}
we just have to show that the sum of the mixed terms \[
(\overline{\partial_{N}}\overline{\partial_{T}}^{*}+\overline{\partial_{T}}^{*}\overline{\partial_{N}})+(\overline{\partial_{N}}^{*}\overline{\partial_{T}}+\partial\overline{_{T}\partial_{N}}^{*})\]
 vanishes. Let us first show that the first term vanishes. Observe
that the following anti-commutation relations hold: \[
d\overline{z_{i}}\wedge\overline{\partial\rho'}^{*}+\overline{\partial\rho'}^{*}d\overline{z_{i}}\wedge=0.\]
Indeed, this is equivalent to the corresponding forms being orthogonal.
Using this and the expansion \ref{II eq exp of dbar} we get that
\[
\overline{(\partial_{N}}\overline{\partial_{T}}^{*}+\overline{\partial_{T}}^{*}\overline{\partial_{N}})=\sum_{i}[\overline{N},\overline{Z}_{i}^{*}]\overline{\partial\rho'}^{*}d\overline{z_{i}}\wedge.\]
But this equals zero since the commutators $[\overline{N},\overline{Z}_{i}^{*}]$
vanish, using the expressions in lemma \ref{II formal adj}. To see
that the second terms vanishes one can go through the same argument
again, now using that the commutators $[\overline{N}^{*},\overline{Z}_{i}]$
vanish. 
\end{proof}
We will call a sequence $\chi_{i}$ of non-negative functions on $X_{0}$
a \emph{relative exhaustion sequence} if there is sequence of balls
$B_{R_{i}}$ centered at the origin and exhausting $\C^{n},$ such
that $\chi_{i}$ is identically equal to 1 on $B_{R_{i}/2}$ and with
support in $B_{R_{i}}.$ Moreover, if the metric $\omega$ is such
that the sequence $\chi_{i}$ can be chosen to make $\left|d\chi_{i}\right|$
uniformly bounded then $(X_{0},\omega)$ is called \emph{relatively
complete}. The point is that when $(X_{0},\omega)$ is relatively
complete, one can integrate partially without getting boundary terms
{}``at infinity''. For a complete manifold this was shown in \cite{gro}
and the extension to the relative case is straightforward. 

\begin{lem}
\label{II lem gromov}Suppose that $(X_{0},\omega)$ is relatively
complete. Then there is a relative exhaustion sequence $\chi_{i}$of
$X_{0}$such that, if $\alpha$is a smooth form in $L^{2}(X_{0}),$then
\[
\lim_{i}(\chi_{i}\Delta_{T}\alpha,\alpha)=(\overline{\partial_{T}}^{*}\alpha,\overline{\partial_{T}}^{*}\alpha)+(\overline{\partial_{T}}\alpha,\overline{\partial_{T}}\alpha).\]
Moreover, if $\alpha$satisfies $\overline{\partial}-$Neumann boundary
conditions on $\partial X_{0}$, then \begin{equation}
\lim_{i}(\chi_{i}\Delta_{N}\alpha,\alpha)=(\overline{\partial_{N}}^{*}\alpha,\overline{\partial_{N}}^{*}\alpha)+(\overline{\partial_{N}}\alpha,\overline{\partial_{N}}\alpha).\label{eq: lemma gromov}\end{equation}

\end{lem}
\begin{proof}
Since $(X_{0},\omega)$ is relatively complete, following section
$1.1B$ in \cite{gro} it is enough to prove the statements for a
form $\alpha$ with relatively compact support, with $\chi_{i}$ identically
equal to $1$ (this is called the Gaffney cutoff trick in \cite{gro}).
Assuming this, the first statement then follows immediately from lemma
\ref{II formal adj}. To prove the second statement we assume that
$\alpha$ satisfies $\overline{\partial}-$Neumann boundary conditions,
i.e. $\alpha_{N}=(\overline{\partial}\alpha)_{N}=0$ on $\partial X_{0}.$
According to lemma \ref{II formal adj} the first term in the right
hand side of \ref{eq: lemma gromov} may be written as $(\overline{\partial_{N}}\overline{\partial_{N}}^{*}\alpha,\alpha),$
since $\alpha_{N}=0$ on $\partial X_{0}$ by assumption. To show
that the second term may be written as $(\overline{\partial_{N}^{*}}\overline{\partial_{N}}\alpha,\alpha)$
we just have to show that $(\overline{\partial_{N}}\alpha)_{N}=0$
on $\partial X_{0}.$ To this end, first observe that $\overline{\partial_{N}}\alpha=\overline{\partial_{N}}\alpha_{T}$
and $(\overline{\partial}\alpha)_{N}=\overline{\partial_{T}}\alpha_{N}+\overline{\partial_{N}}\alpha_{T}.$
Now, by assumption $(\overline{\partial}\alpha)_{N}=0$ on $\partial X_{0}.$
Combing this with the previous two identities we deduce that $\overline{\partial_{N}}\alpha=-\overline{\partial_{T}}\alpha_{N}$
on $\partial X_{0}.$ But $\alpha_{N}=0$ on $\partial X_{0}$ and
$\overline{\partial_{T}}$ is a tangential operator, so it follows
that $\overline{\partial_{T}}\alpha_{N}=0$ on $\partial X_{0}.$
This proves that $(\overline{\partial_{N}}\alpha)_{N}=0$ on $\partial X_{0}.$
\end{proof}
Finally, to finish the proof of proposition \ref{II prop tangential},
first observe that the model metric $\omega_{0}$ corresponding to
$a(\rho_{0})=(1-\rho_{0})^{2}$ is relatively complete. Now take an
arbitrary form $\alpha$ in $\mathcal{H}^{0,q}(X_{0},\phi_{0}).$
Then $(\chi_{i}\Delta\alpha,\alpha)=0$ for each $i.$ Hence, using
lemma \ref{II lem decomp of lapl} together with lemma \ref{II lem gromov}
we deduce, after letting $i$ tend to infinity, that\begin{equation}
0=\left\Vert \overline{\partial_{T}}\alpha\right\Vert ^{2}+\left\Vert \overline{\partial_{T}}^{*}\alpha\right\Vert +\left\Vert \overline{\partial_{N}}\alpha\right\Vert ^{2}+\left\Vert \overline{\partial_{N}}^{*}\alpha\right\Vert ^{2}.\label{eq: sum of squares}\end{equation}
In particular, $\overline{\partial_{N}}^{*}\alpha$ vanishes in $X_{0}.$
If we write $\alpha_{N}=\overline{\partial\rho}'\wedge\sum_{I}g_{I}d\overline{z_{I}}$
this means that \[
\overline{N}^{*}g_{I}=0\]
in $X_{0}$ for all $I.$ Moreover, since $\alpha$ satisfies $\overline{\partial}-$Neumann
boundary conditions, each function $g_{I}$vanishes on the boundary
of $X_{0}.$ It follows that $g_{I}=0$ in all of $X_{0}.$ Indeed,
let $g_{I}':=a(\rho)^{-1/2}\overline{g_{I}}$ and consider the restriction
of $g_{I}'$ to the half planes in $\C$ obtained by freezing the
$z_{i}-$ variables. Then $g_{I}'$ is holomorphic in the half plane,
vanishing on the boundary. It is a classical fact that $g_{I}'$ then
actually vanishes identically. Moreover, \ref{eq: sum of squares}
also gives that $\alpha$ is $\overline{\partial}-$ closed and coclosed.
This finishes the proof of proposition \ref{II prop tangential}.

\section{Contributions from the three boundary regions}

In this section we will estimate the integral of the Bergman function
over the three different boundary regions. The contribution from the
inner part of $X$ was essentially computed in \cite{berm}.

\subsection{\label{II sec first regt}The first region}

Recall that the first region is the set where $\rho\geq-R/k).$ Fix
a point $\sigma$ in $\partial X$ and take local holomorphic coordinates
$(z,w),$ where $z$ is in $\C^{n-1}$ and $w=u+iv.$ By an appropriate
choice, we may assume that the coordinates are orthonormal at $0$
and that \begin{equation}
\rho(z,w)=\sum_{i=1}^{n-1}v-\mu_{i}\left|z_{i}\right|^{2}+O(\left|(z,w)\right|^{3})=:\rho_{0}(z,w)+O(\left|(z,w)\right|^{3}).\label{II local expr for def funct}\end{equation}
 In a suitable local holomorphic trivialization of $L$ close to the
boundary point $\sigma,$ the fiber metric may be written as \[
\phi(z)=\sum_{i,j=1}^{n-1}\lambda_{ij}z_{i}\overline{z_{j}}+O(\left|w\right|)O(\left|z\right|)+O(\left|w\right|^{2})+O(\left|(z,w)\right|^{3}).\]
 Denote by $F_{k}$ the holomorphic scaling map \[
F_{k}(z,w)=(z/k^{1/2},w/k),\]
 so that \[
X_{k}=F_{k}(D_{\ln k})\bigcap X\]
 is a sequence of decreasing neighborhoods of the boundary point $\sigma,$
where $D_{\ln k}$ denotes the polydisc of radius $\ln k$ in $\C^{n}.$
Note that \[
F_{k}^{-1}(X_{k})\rightarrow X_{0},\]
 in a certain sense, where $X_{0}$ is the model domain with defining
function $\rho_{0}.$ On $F_{k}^{-1}(X_{k})$ we have the \emph{scaled
metrics} $F_{k}^{*}k\omega_{k}$ and $F_{k}^{*}k\phi$ that tend to
the model metrics $\omega_{0}$ and $\phi_{0}$ on $X_{0},$ when
$k$ tends to infinity, where \begin{equation}
\omega_{0}=\frac{i}{2}\partial\overline{\partial}\left|z\right|^{2}+a(\rho_{0})^{-1}2i\partial\rho_{0}\wedge\overline{\partial\rho_{0}}\,\,\,\textrm{and}\,\,\,\phi_{0}(z)=\sum_{i,j=1}^{n-1}\lambda_{ij}z_{i}\overline{z_{j}},\label{II model metrics}\end{equation}
 for a smooth function $a(\rho_{0})$ that is positive on $]-\infty,0].$
The factor $a(\rho_{0}),$ and hence the model metric $\omega_{0},$
really depends on the number $R$ (used in the definition of the boundary
regions). However, the dependence on $R$ will play no role in the
proofs, since $R$ will be fixed when $k$ tends to infinity. The
space of {}``model'' harmonic forms in $L^{2}(X_{0},\omega_{0},\phi_{0})$
satisfying $\overline{\partial}-$Neumann boundary conditions will
be denoted by $\mathcal{H}(X_{0},\phi_{0})).$ 

\begin{lem}
\label{II lemma conv of metrics}For the component-wise uniform norms
on $F_{k}^{-1}(X_{k})$ we have that \[
\begin{array}{rcl}
\begin{array}{rcl}
\left\Vert F_{k}^{*}k\rho-\rho_{0}\right\Vert _{\infty} & \rightarrow & 0\\
\left\Vert F_{k}^{*}k\omega_{k}-\omega_{0}\right\Vert _{\infty} & \rightarrow & 0\\
\left\Vert F_{k}^{*}k\phi-\phi_{0}\right\Vert _{\infty} & \rightarrow & 0\end{array}\end{array}\]
and similarly for all derivatives.
\end{lem}
\begin{proof}
The convergence for $\rho$ and $\phi$ is straightforward (compare
\cite{berm}) and the convergence for $\omega_{k}$ will be showed
in section \ref{II subsection metricsequence} once $\omega_{k}$
has been constructed. 
\end{proof}
The Laplacian on $F_{k}^{-1}(X)$ taken with respect to the scaled
metrics will be denoted by $\Delta^{(k)}$ and the corresponding formal
adjoint of $\overline{\partial}$ will be denoted by $\overline{\partial}^{*(k)}.$
The Laplacian on $X_{0}$ taken with respect to the model metrics
will be denoted by $\Delta_{0}.$ Because of the convergence property
of the metrics above it is not hard to check that \begin{equation}
\Delta^{(k)}=\Delta_{0}+\varepsilon_{k}\mathcal{D}_{k},\label{II conv of Laplace}\end{equation}
where $\mathcal{D}_{k}$ is a second order partial differential operator
with bounded variable coefficients on $F_{k}^{-1}(X_{k})$ and $\epsilon_{k}$
is a sequence tending to zero with $k.$ Next, given a $(0,q)-$form
$\alpha_{k}$ on $X_{k}$ with values in $L^{k},$ define the \emph{scaled
form} $\alpha^{(k)}$ on $F_{k}^{-1}(X_{k})$ by \[
\alpha^{(k)}:=F_{k}^{*}\alpha_{k}.\]
Then \begin{equation}
F_{k}^{*}\left|\alpha_{k}\right|^{2}=\left|\alpha^{(k)}\right|^{2},\label{II section first norm scaling}\end{equation}
 where the norm of $\alpha_{k}$ is the one induced by the metrics
$k\omega_{k}$ and $k\phi$ and the norm of the scaled form $\alpha^{(k)}$
is taken with respect to the scaled metrics $F_{k}^{*}k\omega_{k}$
and $F_{k}^{*}k\phi.$ This is a direct consequence of the definitions.
Moreover, the next lemma gives the transformation of the Laplacian.

\begin{lem}
\label{lem: scaling of laplace}The following relation between the
Laplacians holds: \begin{equation}
\Delta^{(k)}\alpha^{(k)}=(\Delta_{k}\alpha_{k})^{(k)}.\label{II eq scaling of laplace}\end{equation}

\end{lem}
\begin{proof}
Since the Laplacian is naturally defined with respect to any given
metric it is invariant under pull-back, proving the lemma. 
\end{proof}
\emph{In the following, all norms over $F_{k}^{-1}(X_{k})$ will be
taken with respect to the model metrics $\omega_{0}$ and $\phi_{0}.$}
The point is that these norms anyway coincide, asymptotically in $k,$
with the norms defined with respect to the scaled metrics used above,
by the following lemma. 

\begin{lem}
\label{II lemma localization of norms etc} We have that uniformly
on $F_{k}^{-1}(X_{k})$ \[
\begin{array}{rcl}
F_{k}^{*}\left|\alpha_{k}\right| & \sim & \left|\alpha^{(k)}\right|\\
\left\Vert \alpha_{k}\right\Vert _{X_{k}} & \sim & \left\Vert \alpha^{(k)}\right\Vert _{F_{k}^{-1}(X_{k})}\end{array}\]
 Moreover,for any sequence $\theta_{k}^{I}$ of $\omega_{k}-$ orthonormal
bases of direction forms in $\Lambda^{0,q}(X)_{x}$ at $F_{k}(x),$
there is a bases of $\omega_{0}-$ orthonormal direction forms at
$x,$ such that the following asymptotic equality holds, when $k$
tends to infinity: \[
F_{k}^{*}\left|\alpha_{k,\theta_{k}^{I}}\right|\sim\left|\alpha_{\theta^{I}}^{(k)}\right|\]
 for each index $I.$ 
\end{lem}
\begin{proof}
The lemma follows immediately from \ref{II section first norm scaling}
and the convergence of the metrics in the previous lemma.
\end{proof}
Now we can prove the following lemma that makes precise the statement
that, in the large $k$ limit, harmonic forms $\alpha_{k}$ are harmonic
with respect to the model metrics and the model domain on a small
scale close to the boundary of $X.$ 

\begin{lem}
\label{II lemma: convergence of forms}Suppose that the boundary of
$X$ satisfies condition $Z(q)$ (see remark \ref{rem: cond z}).
For each $k,$ suppose that $\alpha^{(k)}$ is a $\overline{\partial}-$closed
smooth form on $F_{k}^{-1}(X_{k})$ such that $\overline{\partial}^{*(k)}\alpha^{(k)}=0$
and that $\alpha^{(k)}$ satisfies $\overline{\partial}-$Neumann
boundary conditions on $F_{k}^{-1}(\partial X).$ Identify $\alpha^{(k)}$
with a form in $L^{2}(\C^{n})$ by extending with zero. Then there
is a constant $C_{R}$ independent of $k$ such that 

\[
\sup_{D_{R}\bigcap F_{k}^{-1}(X_{k}).}\left|\alpha^{(k)}\right|^{2}\leq C_{R}\left\Vert \alpha^{(k)}\right\Vert _{D_{2R}\bigcap F_{k}^{-1}(X_{k})}^{2}\]
Moreover, if the sequence of norms $\left\Vert \alpha^{(k)}\right\Vert _{F_{k}^{-1}(X_{k})}^{2}$
is bounded, then there is a subsequence of $\left\{ \alpha^{(k)}\right\} $
which converges uniformly with all derivatives on any compactly included
set in $X_{0}$ to a smooth form $\beta,$ where $\beta$ is in $\mathcal{H}^{0,q}(X_{0}).$
The convergence is uniform on $D_{R}\bigcap F_{k}^{-1}(X_{k}).$
\end{lem}
\begin{proof}
Fix a $k$ and consider the intersection of the polydisc $D_{R}$
of radius $R$ with $F_{k}^{-1}(X_{k}).$ It is well-known that the
Laplace operator $\Delta^{(k)}$ acting on $(0,q)-$ forms is sub-elliptic
close to a point $x$ in the boundary satisfying the condition $Z(q)$
(see \cite{f-k}). In particular, sub-elliptic estimates give for
any smooth form $\beta^{(k)}$ satisfying $\overline{\partial}-$Neumann
boundary conditions on $F_{k}^{-1}(\partial X)$ that \begin{equation}
\left\Vert \beta^{(k)}\right\Vert _{D_{R}\bigcap F_{k}^{-1}(X),m-1}^{2}\leq C_{k,R}(\left\Vert \beta^{(k)}\right\Vert _{D_{2R}\bigcap F_{k}^{-1}(X)}^{2}+\left\Vert \Delta^{(k)}\beta^{(k)}\right\Vert _{D_{2R}\bigcap F_{k}^{-1}(X),m}^{2}),\label{I pf of lemma conv of forms: subelliptic estimate}\end{equation}
 where the subscript $m$ indicates a Sobolev norm with $m$ derivatives
in $L^{2}$ and where the norms are taken in $F_{k}^{-1}(X)$ with
respect to the scaled metrics. The $k-$dependence of the constants
$C_{k,R}$ comes from the boundary $F_{k}^{-1}(\partial X)$ and the
scaled metrics $F_{k}^{*}k\omega_{k}$ and $F_{k}^{*}k\phi.$ However,
thanks to the convergence of the metrics in lemma \ref{II lemma conv of metrics}
one can check that the dependence is uniform in $k.$ Hence, applying
the subelliptic estimates \ref{I pf of lemma conv of forms: subelliptic estimate}
to $\alpha^{(k)}$ we get\begin{equation}
\left\Vert \alpha^{(k)}\right\Vert _{D_{R}\bigcap F_{k}^{-1}(X),m}^{2}\leq C_{R}\left\Vert \alpha^{(k)}\right\Vert _{D_{2R}\bigcap F_{k}^{-1}(X)}^{2}\label{II pf of lemma conv of forms: Sobolev estimate on Fk}\end{equation}
and the continuous injection $L^{2,l}\hookrightarrow C^{0},\, l>n,$
provided by the Sobolev embedding theorem, proves the first statement
in the lemma. To prove the second statement assume that $\left\Vert \alpha^{(k)}\right\Vert _{F_{k}^{-1}(X)}^{2}$
is uniformly bounded in $k.$ Take a sequence of sets $K_{n},$ compactly
included in $X_{0},$ exhausting $X_{0}$ when $n$ tends to infinity.
Then the estimate \ref{I pf of lemma conv of forms: subelliptic estimate}
(applied to polydiscs of increasing radii) shows that \begin{equation}
\left\Vert \alpha^{(k)}\right\Vert _{K_{n},m}^{2}\leq C'_{n}\label{II pf of lemma conv of forms}\end{equation}
 Since this holds for any $m\geq1,$ Rellich's compactness theorem
yields, for each $n,$ a subsequence of $\left\{ \alpha^{(k)}\right\} $,
which converges in all Sobolev spaces $L^{2,l}(K_{n})$ for $l\geq0$
for a fixed $n.$ The compact embedding $L^{2,l}\hookrightarrow C^{p},\, k>n+\frac{1}{2}p,$
shows that the sequence converges in all $C^{p}(K_{n}).$ Choosing
a diagonal sequence with respect to $k$ and $n,$ yields convergence
on any compactly included set $K.$ Finally, we will prove that the
limit form $\beta$ is in $\mathcal{H}(X_{0}).$ First observe that
by weak compactness we may assume that the sequence $1_{X_{0}}\alpha^{(k)}$
tends to $\beta$ weakly in $L^{2}(\C^{n}),$ where $1_{X_{0}}$ is
the characteristic function of $X_{0}$ and $\beta$ is extended by
zero to all of $\C^{n}.$ In particular, the form $\beta$ is weakly
$\overline{\partial}$-closed in $X_{0}.$ To prove that $\beta$
is in $\mathcal{H}(X_{0})$ it will now be enough to show that \begin{equation}
(\beta,\overline{\partial}\eta)_{X_{0}}=0\label{II pf lemma conv of forms coclosed}\end{equation}
for any form $\eta$ in $X_{0}$ that is smooth up to the boundary
and with a relatively compact support in $X_{0}.$ Indeed, it is well-known
that $\beta$ then is in the kernel of the Hilbert adjoint of the
densely defined operator $\overline{\partial.}$ Moreover, the regularity
theory then shows that $\beta$ is smooth up to the boundary, where
it satisfies $\overline{\partial}-$Neumann boundary conditions (actually
this is shown using sub-elliptic estimates as in \ref{I pf of lemma conv of forms: subelliptic estimate})
\cite{h,f-k}. To see that \ref{II pf lemma conv of forms coclosed}
holds, we write the left hand side, using the weak convergence of
$1_{X_{0}}\alpha^{(k)},$ as \begin{equation}
\lim_{k}(\alpha^{(k)},\overline{\partial}\eta)_{X_{0}}=\lim_{k}(\alpha^{(k)},\overline{\partial}\eta)_{X_{0}\bigcap F_{k}^{-1}(X)}.\label{II pf of lemma conv of forms coclose b}\end{equation}
Extending $\eta$ to a smooth form on some neighborhood of $X_{0}$
in $\C^{n}$ we may now write this as a scalar product, with respect
to the scaled metrics, over $F_{k}^{-1}(X),$ thanks to the convergence
in lemma \ref{II lemma conv of metrics} of the scaled metrics and
the scaled defining function. Since, $\alpha^{(k)}$ is assumed to
satisfy $\overline{\partial}-$Neumann-boundary conditions on $F_{k}^{-1}(\partial X)$
and be in the kernel of the formal adjoint of $\overline{\partial},$
taken with respect to the scaled metrics, this means that the right
hand side of \ref{II pf of lemma conv of forms coclose b} vanishes.
This proves \ref{II pf lemma conv of forms coclosed} and finishes
the proof of the lemma. 
\end{proof}
The following proposition will give the boundary contribution to the
holomorphic Morse inequalities in theorem \ref{II thm morse}.

\begin{prop}
\label{I prop first region}Let\[
I_{R}:=\limsup_{k}\int_{-\rho<Rk^{-1}}B_{X}^{q,k}\omega_{n}{}\]
 Then \[
\limsup_{R}I_{R}\leq(-1)^{q}(\frac{1}{2\pi})^{n}\int_{\partial X}\int_{T(q)_{\rho,x}}(\Theta+t\mathcal{L})_{n-1}\wedge\partial\rho\wedge dt\]
 where $T(q)_{\rho,x}=\left\{ t>0:\,\textrm{index$($}\Theta+t\mathcal{L})=q\,\textrm{along }T^{1,0}(\partial X)_{x}\right\} .$
\end{prop}
\begin{proof}
We may assume that the boundary of $X$ satisfies condition $Z(q)$
(compare remark \ref{rem: cond z}). Using the expression \ref{II metric express}
for the metric $\omega_{k},$ the volume form $(\omega_{k})_{n}$
may be written as $a_{k}(\rho)^{-1}(\omega_{T})_{n-1}\wedge2i\partial\rho\wedge d\rho.$
Hence, $I_{R}$ can be expressed as \[
\limsup_{k}\int_{\partial X}(\omega_{T})_{n-1}\wedge2i\partial\rho\int_{-R/k}^{0}a_{k}(\rho)^{-1}B_{X}^{k}d\rho\]
Now, fix a point in the boundary of $X$ and take local coordinates
as in the beginning of the section. To make the argument cleaner we
will first assume that the restriction of $\rho$ to the ray close
to the boundary where $z$ and the real part of $w$ vanish, coincides
with $v.$ Then, after a change of variables, the inner integral along
the ray in the first region becomes \begin{equation}
1/k\int_{0}^{R}a_{k}(v/k)^{-1}B_{X}^{q,k}(v/k)dv.\label{II pf of prop first region integral}\end{equation}
Moreover, by the scaling properties of the metrics $k\omega_{k}$
in lemma \ref{II lemma conv of metrics} we have the uniform convergence
\[
ka_{k}(v/k)\rightarrow a(v).\]
on the segment $[0,R]$ (see also section \ref{II subsection metricsequence}).
Thus, \[
I_{R}=\limsup_{k}\int_{\partial X}(\omega_{T})_{n-1}\wedge2i\partial\rho\int_{0}^{R}a(v)^{-1}B_{X}^{q,k}(v/k)dv.\]
 Let us now show that\begin{equation}
\begin{array}{lrcl}
(i) & B_{X}^{q,k}(0,iv/k) & \lesssim & B_{X_{0}}^{q}(0,iv)\\
(ii) & B_{X}^{q,k}(0,iv/k) & \leq & C_{R}\end{array}\label{II pf of first reg claim}\end{equation}
 We first prove (i). According to the extremal property \ref{II extremal prop of b}
it is enough to show that \begin{equation}
S_{X,\theta_{k}}^{q,k}(0,iv/k)\lesssim S_{X_{0},\theta}^{q}((0,iv)\label{II pr of prop first region claim two}\end{equation}
 for any sequence of direction forms $\theta_{k}$ at $(0,iv/k)$
as in lemma \ref{II lemma localization of norms etc}. Given this,
the bound \ref{II pf of first reg claim} is obtained after summing
over the base elements $\theta_{k}.$ To prove \ref{II pr of prop first region claim two}
we have to estimate 

\[
\left|\alpha_{k,\theta_{k}}(0,iv/k\right|^{2},\]
 where $\alpha_{k}$ is a normalized harmonic form with values in
$L^{k}$ that is extremal at $(0,iv/k)$ in the direction $\theta_{k}.$
Moreover, it is clearly enough to estimate some subsequence of $\alpha_{k}.$
By lemma \ref{II lemma localization of norms etc} it is equivalent
to estimate \[
\left|\alpha_{\theta}^{(k)}(0,iv)\right|^{2},\]
 where the scaled form $\alpha^{(k)}$ is defined on $F_{k}^{-1}(X_{k})$
and extended by zero to all of $X_{0}.$ Note that, according to lemma
\ref{II lemma localization of norms etc} the norms of the sequence
of scaled forms $\alpha^{(k)}$ are asymptotically less than one:
\begin{equation}
\left\Vert \alpha^{(k)}\right\Vert _{F_{k}^{-1}(X)}^{2}\sim\left\Vert \alpha_{k}\right\Vert _{X_{k}}^{2}\leq1,\label{II proof of prop first norm bound}\end{equation}
 since the global norm of $\alpha_{k}$ is equal to one. Hence, by
lemma \ref{II lemma: convergence of forms} there is a subsequence
$\alpha^{(k_{j})}$ that converges uniformly to $\beta$ with all
derivatives on the segment $0\leq v\leq R$ in $X_{0}$ and where
the limit form $\beta$ is in $\mathcal{H}^{0,q}(X_{0},\phi_{0})$
and its norm is less than one (by \ref{II proof of prop first norm bound}).
This means that \[
\left|\alpha_{k,\theta}(0,iv/k)\right|^{2}\sim\left|\alpha_{\theta}^{(k_{j})}(0,iv)\right|^{2}\sim\left|\beta_{\theta}(0,iv)\right|^{2},\]
 Since the limit form $\beta$ is a contender for the model extremal
function $S_{X_{0},\theta}^{q}(0),$ this proves \ref{II pr of prop first region claim two},
and hence we obtain $(i).$ To show $(ii),$ just observe that lemma
\ref{II lemma: convergence of forms} says that there is a constant
$C_{R}$ such that there is a uniform estimate \[
\left|\alpha^{(k)}(0,iv)\right|^{2}\leq C_{R},\]
By the extremal characterization \ref{II extremal prop of b} of $B_{X}^{q,k}$
this proves (ii). Now using \ref{II pf of prop second region claim}
and Fatou's lemma to interchange the limits, $I_{R}$ may be estimated
by \begin{equation}
\int_{\partial X}(\int_{0}^{\infty}B_{X_{0}}^{q}(0,iv)(\omega_{0})_{n},\label{II proof of prop first region integral b}\end{equation}
 in terms of the model metric $\omega_{0}$ on $X_{0}.$ By theorem
\ref{II thm model} this equals \[
(\frac{i}{2\pi})^{n}(-1)^{q}\int_{\partial X}\int_{T(q)}(\partial\overline{\partial}\phi_{0}+t\partial\overline{\partial}\rho_{0})_{n-1}\wedge\partial\rho\wedge dt.\]
This finishes the proof of the proposition under the simplifying assumption
that the restriction of $\rho$ to the ray close to the boundary coincides
with the restriction of $v.$ In general this is only true up to terms
of order $O\left|(z,w)\right|^{3},$ given the expression \ref{II local expr for def funct}.
To handle the general case one writes the integral \ref{II pf of prop first region integral}
as \begin{equation}
1/k\int_{I_{k}}a_{k}B_{X}^{q,k}dv,\label{II proof of prof first region int c}\end{equation}
where $I_{k}$ is the inverse image under $F_{k}$ of the ray in the
first region determined by the fixed point in the boundary. Clearly,
$I_{k}$ tends to the segment $[0,R]$ in $X_{0}$ obtained by keeping
all variables except $v$ equal to zero. Moreover, since the sequence
$\alpha^{(k_{j})}$ above converges uniformly with all derivatives
on $D_{R}\bigcap F_{k}^{-1}(X_{k})$ it forms an equicontinous family,
so the same argument as above gives that \ref{II proof of prof first region int c}
may be estimated by \ref{II proof of prop first region integral b}.
This finishes the proof of the proposition.
\end{proof}

\subsection{The second and third region}

Let us first consider the second region, i.e. where $-1/k^{1/2}\leq\rho\leq-R/k.$
Given a $k,$ consider a fixed point $(0,iv)=(0,ik^{-s}),$ where
$1/2\leq s<1.$ Any point in the second region may be written in this
way. Let $(z',w')$ be coordinates on the unit polydisc $D.$ Define
the following holomorphic map from the unit polydisc $D$ to a neighborhood
of the fixed point: \[
F_{k,s}(z',w')=(k^{-1/2}z',k^{-s}+\frac{1}{2}k^{-s}w')\]
so that \[
X_{k,s}:=F_{k,s}(D),\]
 is a neighborhood of the fixed point, staying away from the boundary
of $X.$ On $D$ we will use the scaled metrics $F_{k,s}^{*}k\omega_{k}$
and $kF_{k,s}^{*}k\phi$ that have bounded derivatives and are comparable
to flat metrics in the following sense: \begin{equation}
\begin{array}{rcccl}
C^{-1}\omega_{E} & \leq & F_{k,s}^{*}k\omega_{k} & \leq & C\omega_{E}\\
 &  & \left|F_{k,s}^{*}k\phi\right| & \leq & C\end{array}\label{II second region equivalence of metrics}\end{equation}
where $\omega_{E}$ is the Euclidean metric. Note that the scaling
property of $\omega_{k}$ is equivalent to \begin{equation}
C^{-1}k\rho^{2}\leq a_{k}(\rho)\leq Ck\rho^{2}.\label{II second region scaling of conf}\end{equation}
 These properties will be verified in section \ref{II subsection metricsequence}
once $\omega_{k}$ is defined. The Laplacian on $D$ taken with respect
to the scaled metrics will be denoted by $\Delta^{(k,s)}$ and the
Laplacian on $X_{0}$ taken with respect to the model metrics will
be denoted by $\Delta_{0}.$ Next, given a $(0,q)-$form $\alpha_{k}$
on $X_{k}$ with values in $L^{k},$ define the scaled form $\alpha^{(k)}$
on $D$ by \[
\alpha^{(k,s)}:=F_{k,s}^{*}\alpha_{k}.\]
Using \ref{II second region equivalence of metrics} one can see that
the following equivalence of norms holds: \begin{equation}
\begin{array}{rcccl}
C^{-1}\left|\alpha^{(k,s)}\right|^{2} & \leq & F_{k,s}^{*}\left|\alpha_{k}\right|^{2} & \leq & C\left|\alpha^{(k,s)}\right|^{2}\\
C^{-1}\left\Vert \alpha^{(k,s)}\right\Vert _{D}^{2} & \leq & \left\Vert \alpha_{k}\right\Vert _{X_{k}}^{2} & \leq & C\left\Vert \alpha^{(k,s)}\right\Vert _{D}^{2}\end{array}\label{II second region equivalence of norms}\end{equation}
 \emph{In the following, all norms over $F_{k,s}^{-1}(X_{k,s})$ will
be taken with respect to the Euclidean metric $\omega_{E}$ and the
trivial fiber metric. }

\begin{prop}
Let \[
II_{R}:=\limsup_{k}\int_{Rk^{-1}<-\rho<k^{-1/2}}B_{X}^{q,k}(\omega_{k})_{n}\]
Then \[
\lim_{R}II_{R}=0.\]

\end{prop}
\begin{proof}
Fix a $k$ and a point $(0,iv)=(0,ik^{-s}),$ where $s$ is in $[1/2,1[.$
From the scaling properties \ref{II second region scaling of conf}
of $\omega_{k}$ it follows that at the point $(0,ik^{-s})$ \[
\omega_{k}^{n}\leq Ck^{(2s-1)}\omega^{n}.\]
 Next, observe that \begin{equation}
B_{X}^{q,k}(0,ik^{-s})\leq C\label{II pf of prop second reg: claim}\end{equation}
 Accepting this for the moment, it follows that \[
B^{q,k}(0,iv)\omega_{k}^{n}\leq Ck^{-1}v^{-2}\omega^{n},\]
since we have assumed that $v=k^{-s}.$ Hence, the integral in $II_{R}$
may be estimated by \[
\int_{\partial X}Ck^{-1}\int_{-k^{-1/2}}^{-Rk^{-1}}v^{-2}\textrm{d}v=Ck^{-1}(R^{-1}k-k^{1/2})\]
 which tends to $CR^{-1}$ when $k$ tends to infinity. This proves
that $II_{R}$ tends to zero when $R$ tends to infinity, which proves
the proposition, given \ref{II pf of prop second reg: claim}. 

Finally, let us prove the claim \ref{II pf of prop second reg: claim}.
For a given $k$ consider the point $(0,ik^{-s})$ as above. As in
the proof of the previous proposition we have to prove the estimate
\begin{equation}
\left|\alpha_{k}(0,ik^{-s})\right|^{2}\leq C,\label{II pf of prop second region claim}\end{equation}
 where $\alpha_{k}$ is a normalized harmonic section with values
in $L^{k}$ that is extremal at $(0,ik^{-s}).$ By the equivalence
of norms \ref{II second region equivalence of norms}, it is equivalent
to prove \[
\left|\alpha^{(k,s)}(0)\right|^{2}\leq C,\]
 where the scaled form $\alpha^{(k,s)}$ is defined on $D.$ Note
that, according to \ref{II second region equivalence of norms} \begin{equation}
\left\Vert \alpha^{(k,s)}\right\Vert _{D}^{2}\sim\left\Vert \alpha_{k}\right\Vert _{X_{k}}^{2}\leq C,\label{II pf of prop second region bound on norm}\end{equation}
 since the global norm of $\alpha_{k}$ is equal to $1.$ Moreover,
a simple modification of lemma \ref{lem: scaling of laplace} gives\[
\Delta^{(k,s)}\alpha^{(k,s)}=0\]
 on $D.$ Since $\Delta^{(k,s)}$ is an elliptic operator on the polydisc
$D,$ inner elliptic estimates (i.e. Gårding's inequality) and the
Sobolev embedding theorem can be used as in \cite{berm} to get \[
\left|\alpha^{(k)}(0)\right|^{2}\leq C\left\Vert \alpha^{(k)}\right\Vert _{D}^{2},\]
 where the constant $C$ is independent of $k$ thanks to the equivalence
\ref{II second region equivalence of metrics} of the metrics. Using
\ref{II pf of prop second region bound on norm}, we obtain the claim
\ref{II pf of prop second region claim}.
\end{proof}
Let us now consider the third region where $-\varepsilon<\rho<-k^{-1/2}.$

\begin{prop}
Let\[
III_{\varepsilon}:=\int_{k^{-1/2}<-\rho<\varepsilon}B_{X}^{q,k}(\omega_{k})_{n}\]
 Then \[
III_{\varepsilon}=O(\varepsilon).\]
 
\end{prop}
\begin{proof}
We just have to observe that\[
k^{-n}B_{X}^{q,k}\leq C,\]
when $\rho<-k^{-1/2}.$ This follows from inner elliptic estimates
as in the proof of the previous proposition, now using $s=1/2$ (compare
\cite{berm}).
\end{proof}

\subsection{End of the proof of theorem \ref{II thm morse} (the weak Morse inequalities)}

First observe that \begin{equation}
\int_{0<-\rho<\varepsilon}B_{X}^{q,k}(\omega_{k})_{n}\lesssim(-1)^{q}(\frac{1}{2\pi})^{n}\int_{\partial X}\int_{T(q)}dt(\partial\overline{\partial}\phi+t\partial\overline{\partial}\rho)^{n-1}\wedge\partial\rho+o(\varepsilon).\label{II pf of Morse first int}\end{equation}
 Indeed, for a fixed $R$ we may write the limit of integrals above
as the sum $I_{R}+II_{R}+III_{\varepsilon}.$ Letting $R$ tend to
infinity and using the previous three propositions we get the estimate
above. Moreover, we have that \begin{equation}
\int_{X_{\varepsilon}}B_{X}^{q,k}\lesssim\int_{X_{\varepsilon}}(\partial\overline{\partial}\phi)^{n},\label{II pf of Morse}\end{equation}
 where $X_{\varepsilon}$ denotes the set where $-\rho$ is larger
than $\epsilon$. This follows from the estimates \[
B_{X}^{q,k}\omega_{n}\lesssim(\frac{i}{2\pi})^{n}(-1)^{q}1_{X(q)}(\partial\overline{\partial}\phi)_{n}\,\,\,\textrm{and}\,\,\, B_{X}^{q,k}\leq\textrm{C}\,\,\,\textrm{in X}_{\varepsilon},\]
 proved in \cite{berm}. Finally, writing $\dim_{\C}\mathcal{H}^{0,q}(X,L^{k})$
as the sum of the integrals in \ref{II pf of Morse first int} and
\ref{II pf of Morse} and letting $\epsilon$ tend to zero, yields
the dimension bound in theorem \ref{II thm morse} for the space of
harmonic forms. By the Hodge theorem we are then done.

\subsection{\label{II subsection metricsequence}The sequence of metrics $\omega_{k}$}

In this section the metrics $\omega_{k}$ will be defined and their
scaling properties, that were used above, will be verified. Recall
that we have to define a sequence of smooth functions $a_{k}$ such
that the metrics \[
\omega_{k}=\omega_{T}+a_{k}(\rho)^{-1}2i\partial\rho\wedge\overline{\partial\rho},\]
 have the scaling properties of lemma \ref{II lemma conv of metrics}
in the first region and satisfy \ref{II second region equivalence of metrics}
in the second region. First observe that the tangential part $\omega_{T}$
clearly scales the right way, i.e. that $F_{k}^{*}k\omega_{T}$ tends
to $\frac{i}{2}\partial\overline{\partial}\left|z\right|^{2}.$ Indeed,
since the coordinates $(z,w)$ are orthonormal at $0$ the forms $\omega_{T}$
and $\frac{i}{2}\partial\overline{\partial}\left|z\right|^{2}$ coincide
at $0.$ Since $\frac{i}{2}\partial\overline{\partial}\left|z\right|^{2}$
is invariant under $F_{k}^{*}k$ the convergence then follows immediately.
We now consider the normal part of $\omega_{k}$ and show how to define
the functions $a_{k}.$ Consider first the piecewise smooth functions
$\widetilde{a_{k}}$ where $\widetilde{a_{k}}$ is defined as $R^{2}/k$
in the first region, as $k\rho^{2},$ in the second region and as
$1$ in the third region and on the rest of $X.$ Then it is not hard
to check that $\widetilde{a_{k}}$ satisfies our demands, except at
the two middle boundaries between the three regions. We will now construct
$a_{k}$ as a regularization of $\widetilde{a_{k}}.$ To this end
we write $\widetilde{a_{k}}=k\widetilde{b_{k}}^{2},$ where $\widetilde{b_{k}}$
is defined by \[
\left\{ \begin{array}{lc}
Rk^{-1}, & -\rho\leq Rk^{-1}\\
-\rho, & Rk^{-1}\leq-\rho\leq k^{-1/2}\\
k^{-1/2}, & -\rho\leq k^{-1/2}\end{array}\right.\]
 in the three regions. It will be enough to regularize the sequence
of continuous piecewise linear functions $\widetilde{b}_{k}.$ Decompose
$\widetilde{b}_{k}$ as a sum of continuous piecewise linear functions
\[
\widetilde{b}_{k}(-\rho)=\frac{R}{k}\widetilde{b}_{1,k}(-\frac{k}{R}\rho)+\frac{1}{k^{1/2}}\widetilde{b}_{2,k}(-k^{1/2}\rho),\]
 where $\widetilde{b}_{1,k}$ is determined by linearly interpolating
between \[
\begin{array}{ccc}
\widetilde{b}_{1,k}(0)=1 & \widetilde{b}_{1,k}(1)=1 & \widetilde{b}_{1,k}(k^{1/2}/R)=0\end{array}\,\widetilde{b}_{1,k}(\infty)=0\]
 and $\widetilde{b}_{2,k}$ is determined by \[
\begin{array}{ccc}
\widetilde{b}_{2,k}(0)=0 & \widetilde{b}_{2,k}(Rk^{-1/2})=0 & \widetilde{b}_{2,k}(1)=1\end{array}\,\widetilde{b}_{2,k}(\infty)=1\]
 Now consider the function $b_{k}$ obtained by replacing $\widetilde{b}_{1,k}$
and $\widetilde{b}_{2,k}$ with the continuous piecewise linear functions
$b_{1}$ and $b_{2},$ where $b_{1}$ is determined by\[
\begin{array}{ccc}
b_{1}(0)=1 & b_{1}(1/2)=1 & b_{1}(1)=0\end{array}\, b_{1}(\infty)=0,\]
 and $b_{2}$ is determined by\[
\begin{array}{ccc}
b_{2}(0)=0 & b_{2}(1)=1 & \, b(\infty)=1\end{array}.\]
 Finally, we smooth the corners of the two functions $b_{1}$ and
$b_{2}.$ Let us now show that the sequence of regularized functions
$b_{k}$ scales in the right way. In the first region we have to prove
that lemma \ref{II lemma conv of metrics} is valid, which is equivalent
to showing that there is a function $b_{0}$ such that \begin{equation}
kb_{k}(t/k)\rightarrow b_{0}\label{II approx first}\end{equation}
 with all derivates, for $t$ such that $0\leq t\leq\textrm{Rln}k.$
From the definition we have that \begin{equation}
kb_{k}(t/k)=Rb_{1}(t/R)+t,\label{II reg of metric b}\end{equation}
which is even independent of $k,$ so \ref{II approx first} is trivial
then. Next, consider the second region. To show that \ref{II second region equivalence of metrics}
holds we have to show that, for parameters $s$ such that $1/2\leq s<1,$
the $t-$dependent functions $k^{s}b_{k}(1/k^{s}+t/2k^{s})$ (where
$\left|t\right|\leq1)$ are uniformly bounded from above and below
by positive constants and have uniformly bounded derivatives. First
observe that in the second region the sequence of functions may be
written as \[
k^{s-1/2}b_{2}(1/k^{1/2-s}+t/2k^{1/2-s}),\]
 and it is not hard to see that it is bounded from above and below
by positive constants, independently of $s$ and $k.$ Moreover, differentiating
with respect to $t$ shows that all derivatives are bounded, independently
of $s$ and $k.$ All in all this means that we have constructed a
sequence of metrics $\omega_{k}$ with the right scaling properties.
In particular, \ref{II reg of metric b} shows that the factor $a^{-1}(\rho_{0})$
in the model metric $\omega_{0}$ \ref{II model metrics} satisfies
\[
C_{R}^{-1}(1-\rho_{0})^{-2}\leq a^{-1}(\rho_{0})\leq C_{R}(1-\rho_{0})^{-2}\]
 for some constant $C_{R}$ depending on $R.$

\part{The strong Morse inequalities and sharp examples}

\section{The strong Morse inequalities}

We will assume that the boundary of $X$ satisfies condition $Z(q)$
(compare remark \ref{rem: cond z}). We will use the same notation
as in section \ref{II section harmonic}. Let $\mu_{k}$ be a sequence
tending to zero. Denote by $\mathcal{H}_{\leq\mu_{k}}^{0,q}(X)$ the
space spanned by the $(0,q)-$eigenforms of the Laplacian $\Delta,$
with eigenvalues bounded by $\mu_{k}.$ The forms are assumed to satisfy
$\overline{\partial}-$Neumann boundary conditions and they will be
called \emph{low energy forms.} Since we have assumed that condition
$Z(q)$ holds, this space is finite dimensional for each $k$ \cite{f-k}.
Recall that the Laplacian is defined with respect to the metric $k\omega_{k},$
so that the eigen values corresponding to $\mu_{k}$ are multiplied
with $k$ if the metric $\omega_{k}$ is used instead.

We will first show that the weak holomorphic Morse inequalities are
\emph{equalities} for the space $\mathcal{H}_{\leq\mu_{k}}^{0,q}(X)$
of low energy forms. When $X$ has no boundary this yields \emph{strong}
Morse inequalities for the truncated Euler characteristics of the
Dolbeault complex with values in $L^{k}.$ However, when $X$ has
a boundary one has to assume that the boundary of $X$ has either
concave or convexity properties to ensure that the corresponding cohomology
groups are finite dimensional, in order to obtain strong Morse inequalities.

The Bergman form for the space $\mathcal{H}_{\leq\mu_{k}}^{0,q}(X)$
defined as in section \ref{II section bergman} will be denoted by
$\B_{\leq\mu_{k}}^{q}.$ By $L_{m}^{2}(X)$ we will denote the Sobolev
space with $m$ derivatives in $L^{2}(X)$ and a subscript $m$ on
a norm will indicate the corresponding Sobolev norm. The essential
part in proving that we now get equality in the weak Morse inequalities
is to show that the estimate on the Bergman form \ref{II pf of first reg claim}
in the proof of proposition \ref{I prop first region} becomes an
asymptotic equality, when considering low energy forms. The rest of
the argument is more or less as before. 

Let us first prove the upper bound, i.e. that the low energy Bergman
form $\B_{\leq\mu_{k}}^{q}$ is asymptotically bounded by the model
harmonic Bergman form. 

\begin{prop}
\label{II prop upper bound on lowenergy b}We have that \[
B_{\leq\mu_{k},\theta_{k}}^{q}(0,iv/k)\lesssim B_{X_{0},\theta}(0,v)\]
 and the sequence $B_{\leq\mu_{k}}^{q}(0,iv/k)$ is uniformly bounded
in the first region.
\end{prop}
\begin{proof}
Let $\alpha_{k}$ be a sequence of normalized forms, such that $\alpha_{k}$
is an extremal for the Hilbert space $\mathcal{H}_{\leq\mu_{k}}^{0,q}(X)$
at the point $(0,iv/k)$ in the direction $\theta_{k}.$ In the following
all norms will be taken over $F_{k}^{-1}(X).$ Observe that by the
invariance property in lemma \ref{lem: scaling of laplace} of the
Laplacian, the scaled form $\alpha^{(k)}$ satisfies\begin{equation}
\left\Vert (\Delta^{(k)})^{p}\alpha^{(k)}\right\Vert ^{2}\leq\mu_{k}^{2p}\rightarrow0\label{II pf of prop upper bd on b}\end{equation}
 for all positive integers $p.$ Let us now show that \begin{equation}
\left\Vert \Delta^{(k)}\alpha^{(k)}\right\Vert _{m}^{2}\rightarrow0\label{II pf of prop upper bound on b b}\end{equation}
 for all non-negative integers $m.$ First observe that $(\Delta^{(k)})^{p}\alpha^{(k)}$
satisfies $\overline{\partial}-$Neumann boundary conditions for all
$p.$ Indeed, by definition all forms in the space $\mathcal{H}_{\leq\mu_{k}}(X_{0})$
satisfy $\overline{\partial}-$Neumann boundary conditions and since
$\Delta$ preserves this space, the forms $(\Delta)^{p}\alpha_{k}$
also satisfy $\overline{\partial}-$Neumann boundary conditions for
all $p.$ By the scaling of the Laplacian this means that the forms
$(\Delta^{(k)})^{p}\alpha^{(k)}$ satisfy $\overline{\partial}-$Neumann
boundary conditions with respect to the scaled metrics. Now applying
the subelliptic estimates \ref{I pf of lemma conv of forms: subelliptic estimate}
to forms of the type $(\Delta^{(k)})^{p}\alpha^{(k)}$ one gets ,
using induction, that \[
\left\Vert \Delta^{(k)}\alpha^{(k)}\right\Vert _{m}^{2}\leq C\sum_{j=1}^{m+1}\left\Vert (\Delta^{(k)})^{j}\alpha^{(k)}\right\Vert ^{2}.\]
 Combining this with \ref{II pf of prop upper bd on b} proves \ref{II pf of prop upper bound on b b}.
Now the rest of the argument proceeds almost word for word as in the
proof of the claim \ref{II pf of first reg claim} in the proof of
proposition \ref{I prop first region}. The point is that the limit
form $\beta$ will still be in $\mathcal{H}(X_{0}),$ thanks to \ref{II pf of prop upper bd on b}. 
\end{proof}
Let us now show how to get the corresponding reverse bound for $\B_{\leq\mu_{k}}^{q,k}$.
We will have use for the following lemma.

\begin{lem}
Suppose that $\beta$ is a normalized extremal form for $\mathcal{H}^{0,q}(X_{0},\phi_{0})$
at the point $(0,iv_{0})$ in the direction $\theta.$ Then \begin{equation}
\left|\beta_{\theta}(0,iv_{0})\right|^{2}=\frac{1}{4\pi}\int_{T(q)}B_{t,\theta}(z,z)e^{v_{0}t}b(t)^{-1}dt.\label{II lemma extremal model statem}\end{equation}
 with notation as in lemma \ref{II lemma model k expansion}. Moreover,
$\beta$ is in $L_{m}^{2}(X_{0})$ for all $m.$ 
\end{lem}
\begin{proof}
Let $x_{0}$ be the point $(0,iv_{0})$ in $X_{0}.$ Since $\beta$
is extremal we have, according to section \ref{II section bergman},
that $\left|\beta_{\theta}(0,iv_{0})\right|^{2}=B_{X_{0},\theta}(0,iv_{0}),$
which in turn gives \ref{II lemma extremal model statem} according
to lemma \ref{II lemma model k expansion}. To prove that $\beta$
is in $L_{m}^{2}(X_{0})$ for all $m,$ we write $\beta$ as \[
\beta(z,w)=\int_{T(q)}\widehat{\beta_{t}}(z)e^{\frac{i}{2}wt}dt,\]
in terms of its Fourier transform as in section \ref{II section model}.
Recall that we have assumed that condition $Z(q)$ holds on the boundary
of $X,$ so that $T(q)$ is finite. Using proposition \ref{II lemma model scalar pr}
we can write \[
\left\Vert \frac{\partial^{l}}{\partial^{l}w}\partial^{I\overline{J}}\beta\right\Vert _{X_{0}}^{2}=4\pi\int_{T(q)}\left\Vert \partial^{I\overline{J}}\widehat{\beta_{t}}(z)\right\Vert _{t}^{2}t^{2l}b(t)dt,\]
where $\partial^{I\overline{J}}$ denotes the complex partial derivatives
taken with respect to $z_{i}$and $\overline{z_{j}}$ for $i$ and
$j$ in the multi index set $I$ and $J,$ respectively. Since, by
assumption, $\beta$ is in $L^{2}(X_{0})$ the integral converges
for $l=0$ with $I$ and $J$ empty. Now it is enough to show that
$\widehat{\beta_{t}}$ is in $L_{m}^{2}(\C^{n},t\psi+\phi)$ for all
$t$ and positive integers $m.$ To this end we will use the following
generalization of \ref{II extremal prop of b direction}:\begin{equation}
\left|\K_{x,\theta}(y)\right|^{2}=\left|\beta(y)\right|^{2}B_{\theta}(x)\label{(II) proof of lemma: global statement}\end{equation}
 if $\beta$ is an extremal form at the point $x,$ in the direction
$\theta$ (compare \cite{berm2}). By lemma \ref{II lemma model k expansion}
the Fourier transform of $\K_{x,\theta}$ evaluated at $t$ is proportional
to $\K_{z,\theta,t}$ where $\K_{z,\theta,t}$ is the Bergman kernel
form for the space $\mathcal{H}^{0,q}(\C^{n-1},t\psi+\phi)$ at the
point $z$ (and $x=(z,w)$) in the direction $\theta$. In \cite{berm}
it was essentially shown that $\K_{z,\theta,t}$ is in $L_{m}^{2}(\C^{n},t\psi+\phi)$
(more precisely: the property was shown to hold for the corresponding
extremal form). Hence, the same thing holds for $\widehat{\beta_{t}},$
according to \ref{(II) proof of lemma: global statement}, which finishes
the proof of the lemma. 
\end{proof}
Now we can construct a sequence $\alpha_{k}$ of approximate extremals
for the space $\mathcal{H}_{\leq\mu_{k}k}(X)$ of low energy forms. 

\begin{lem}
\label{II lemma approx extremal}For any point $x_{0,k}=(0,iv_{0}/k)$
and direction form $\theta$ in the first region there is a sequence
$\left\{ \alpha_{k}\right\} $ and direction forms $\theta_{k}$such
that $\alpha_{k}$ is in $\Omega^{0,q}(X,L^{k})$ and \[
\begin{array}{lrcl}
(i) & \left|\alpha_{k,\theta_{k}}(0,iv_{0}/k)\right|^{2} & \sim & B_{X_{0},\theta}(0,iv_{0})\\
(ii) & \left\Vert \alpha_{k}\right\Vert _{X}^{2} & \sim & 1\\
(iii) & \left\Vert (\overline{\partial}+\overline{\partial}^{*(k)})\alpha^{(k)}\right\Vert _{m}^{2} & \sim & 0\\
(iv) & (\Delta\alpha_{k},\alpha_{k})_{X} & \leq & \delta_{k}\left\Vert \alpha_{k}\right\Vert ^{2}\end{array}\]
where $\delta_{k}$ is a sequence, independent of $x_{0,k},$ tending
to zero, when $k$ tends to infinity. Moreover, $\alpha_{k}$ satisfies
$\overline{\partial}-$Neumann boundary conditions on $\partial X\bigcap F_{k}^{-1}(D),$
where $D$ is a polydisc in $\C^{n}$centered at $0.$
\end{lem}
\begin{proof}
Consider a sequence of points $x_{0,k}$ that can be written as $(0,iv_{0}/k)$
in local coordinates as in section \ref{II sec first regt}. Let us
first construct a form $\alpha_{k}$ with the properties $(i)$ to
$(iv).$ It is defined by \[
\alpha_{k}:=(F_{k}^{-1})^{*}(\chi_{k}\beta)\]
 where $\chi_{k}(z,w)=\chi(z/\ln k,w/\ln k)$ for $\chi$ a smooth
function in $\C^{n}$ that is equal to one on the polydisc $D$ of
radius one centered at $0,$ vanishing outside the polydisc of radius
two and where $\beta$ is the extremal form at the point $(0,v_{0})$
in the direction $\theta$ from the previous lemma. The definition
of $\alpha_{k}$ is made so that\[
\alpha^{(k)}=\chi_{k}\beta\]
We have used the fact that the form $\beta$ extends naturally as
a smooth form to the domain $X_{\delta}$ with defining function $\rho_{0}-\delta,$
to make sure that $\alpha_{k}$ is defined on all of $X.$ The extension
is obtained by writing $\beta$ in terms of its Fourier transform
with respect to $t$ as in the proof of the previous lemma: \[
\beta(z,w)=\int_{T(q)}\widehat{\beta_{t}}(z)e^{-\frac{i}{2}wt}dt,\]
In fact, the right hand side is defined for all $w$ since we have
assumed that condition $Z(q)$ holds on the boundary so that $T(q)$
is finite. Note that the $L_{m}^{2}-$norm of $\beta$ over $X_{\delta}$
tends to the $L_{m}^{2}-$norm of $\beta$ over $X_{0}$ when $\delta$
tends to zero, as can be seen from the analog of proposition \ref{II lemma model scalar pr}
on the domain $X_{\delta}.$ Indeed, the dependence on $\delta$ only
appears in the definition of $b(t),$ where the upper integration
limit is shifted to $\delta.$ Now the statements $(i)$ and $(ii)$
follow from the corresponding statements in the previous lemma. To
see that $(iii)$ holds, first observe that \[
\overline{\partial}^{*(k)}=\overline{\partial}^{*0}+\varepsilon_{k}\mathcal{D},\]
 where $\mathcal{D}$ \emph{}is a first order differential operator
with bounded coefficients on the ball \emph{$B_{lnk}(0)$} and \emph{}$\varepsilon_{k}$
is a sequence tending to zero. Indeed, this is a simple modification
of the statement \ref{II conv of Laplace}. Moreover, by construction
$(\overline{\partial}+\overline{\partial}^{*0})\beta=0.$ Hence, Leibniz
rule gives \[
\left\Vert (\overline{\partial}+\overline{\partial}^{*(k)})\alpha^{(k)}\right\Vert \leq\delta_{k}\left\Vert \beta\right\Vert _{1}+\left\Vert \left|d\chi_{k}\right|\beta\right\Vert \]
 The first term tends to zero since $\beta$ is in $L_{1}^{2}(X_{0})$
and the second terms tends to zero, since it can be dominated by the
{}``tail'' of a convergent integral. The estimates for $m\geq1$
are proved in a similar way (compare \cite{berm}). Finally, to prove
$(iv)$ observe that by the scaling property \ref{II eq scaling of laplace}
for the Laplacian \[
(\Delta\alpha_{k},\alpha_{k})_{X}=\left\Vert (\overline{\partial}+\overline{\partial}^{*})\alpha_{k}\right\Vert _{X}^{2}=\left\Vert (\overline{\partial}+\overline{\partial}^{*(k)})\alpha^{(k)}\right\Vert .\]
 By $(ii),$ the norm of $\left\Vert \alpha_{k}\right\Vert _{X}^{2}$
tends to one and the norm in the right hand side above can be estimated
as above. To see that $\delta_{k}$ can be taken to be independent
of the point $x_{0,k}$ in the first region, one just observes that
the constants in the estimates depend continuously on the eigenvalues
of the curvature forms (compare \cite{berm}). Finally, consider a
polydisc $D$ in $\C^{n}$ with small radius. We will perturb $\alpha_{k}$
slightly so that it satisfies $\overline{\partial}-$Neumann boundary
conditions on $\partial X\bigcap F_{k}^{-1}(D)$ while preserving
the properties $(i)$ to $(iv).$ Recall that a form $\overline{\partial}-$closed
form $\eta_{k}$ satisfies $\overline{\partial}-$Neumann boundary
conditions on $\partial X$ if \begin{equation}
\overline{\partial\rho}^{*}\eta_{k}=0,\label{II eq neumann}\end{equation}
where $\overline{\partial\rho}^{*}$ is the fiber-wise adjoint of
the operator obtained by wedging with the form $\overline{\partial\rho},$
and where the adjoint is taken with respect to the metric $\omega_{k}$
on $X.$ Equivalently, \[
\overline{\partial k\rho^{(k)}}^{*}\eta^{(k)}=0,\]
 where the adjoint is taken with respect to the scaled metrics. By
construction we have that $\alpha^{(k)}$ is $\overline{\partial}-$closed
on $F_{k}^{-1}(D)$ and \begin{equation}
\overline{\partial\rho_{0}}^{*0}\alpha^{(k)}=0,\label{II eq neumann model}\end{equation}
 where now the adjoint is taken with respect to the model metrics.
Let \[
u^{(k)}:=-\overline{\partial}(k\rho^{(k)}\overline{\partial k\rho^{(k)}}^{*}\alpha^{(k)})\]
and let $\widetilde{\alpha}^{(k)}:=\alpha^{(k)}+\chi u^{(k)},$ where
$\chi$ is the cut-off function defined above. Then, using that $\rho$
vanishes on $\partial X,$ we get that the $\overline{\partial}-$closed
form $\widetilde{\alpha}^{(k)}$ satisfies the scaled $\overline{\partial}-$Neumann
boundary conditions, i.e. the relation \ref{II eq neumann} on $\partial X.$
Moreover, using that $k\rho^{(k)}$ converges to $\rho$ with all
derivatives on a fixed polydisc centered at $0$ (lemma \ref{II lemma conv of metrics})
and \ref{II eq neumann model} one can check that $u^{(k)}$ tends
to zero with all derivatives in $X_{\delta}.$ Finally, since $\chi u^{(k)}$
is supported on a bounded set in $X_{\delta}$ and converges to zero
with all derivatives it is not hard to see that $\widetilde{\alpha_{k}}$
also satisfies the properties $(i)$ to $(iv),$ where $\widetilde{\alpha_{k}}:=F_{k}^{-1}{}^{*}(\widetilde{\alpha^{(k)}})$ 
\end{proof}
By projecting the sequence $\alpha_{k}$ of approximate extremals,
from the previous lemma, on the space of low energy forms we will
now obtain the following lower bound on $\B_{\leq\mu_{k}}^{q}.$

\begin{prop}
\label{II prop lower bound on b low}There is a sequence $\mu_{k}$
tending to zero such that \[
\liminf_{k}B_{\leq\mu_{k},\theta_{k}}^{q}(0,iv/k)\geq B_{X_{0}}(0,v)_{\theta}.\]

\end{prop}
\begin{proof}
The proof is a simple modification of the proof of proposition 5.3
in \cite{berm}. Let $\left\{ \alpha_{k}\right\} $ be the sequence
that the previous lemma provides and decompose it with respect to
the orthogonal decomposition $\Omega^{0,q}(X,L^{k})=\mathcal{H}_{\leq\mu_{k}}^{q}(X,L^{k})\oplus\mathcal{H}_{>\mu_{k}}^{q}(X,L^{k}),$
induced by the spectral decomposition of the subelliptic operator
$\Delta$ \cite{f-k}: \[
\alpha_{k}=\alpha_{1,k}+\alpha_{2,k}\]
First, we prove that \begin{equation}
\lim_{k}\left|\alpha_{2}^{(k)}(0,iv)\right|^{2}=0.\label{II pf of prop lower bound on lowe claim}\end{equation}
Since $\alpha_{2}^{(k)}=\alpha^{(k)}-\alpha_{1}^{(k)}$ the form $\alpha_{2}^{(k)}$
satisfies $\overline{\partial}-$Neumann boundary conditions on the
intersection of the polydisc $D$ with $F_{k}^{-1}(\partial X),$
using lemma \ref{II lemma approx extremal}. Subelliptic estimates
as in the proof of lemma \ref{II prop upper bound on lowenergy b}
then show that\begin{equation}
\left|\alpha_{2}^{(k)}(0)\right|^{2}\leq C(\left(\left\Vert \alpha_{2}^{(k)}\right\Vert _{D\bigcap F_{k}^{-1}(X)}^{2}+\left\Vert (\Delta^{(k)})\alpha_{2}^{(k)}\right\Vert _{D\bigcap F_{k}^{-1}(\partial X),m}^{2}\right)\label{II pf of prop upper bound on lowe b}\end{equation}
 for some large integer $m.$ To see that the first term in the right
hand side tends to zero, we first estimate $\left\Vert \alpha_{2}^{(k)}\right\Vert _{D\bigcap F_{k}^{-1}(X)}^{2}$
with $\left\Vert \alpha_{k,2}\right\Vert _{X}^{2}$ using the norm
localization in lemma \ref{II lemma localization of norms etc}. Next,
by the spectral decomposition of $\Delta_{k}:$\[
\left\Vert \alpha_{2,k}\right\Vert _{X}^{2}\leq\frac{1}{\mu_{k}}\left\langle \Delta_{k}\alpha_{2,k},\alpha_{2,k}\right\rangle _{X}\leq\frac{1}{\mu_{k}}\left\langle \Delta_{k}\alpha_{k},\alpha_{k}\right\rangle _{X}\leq\frac{\delta_{k}}{\mu_{k}}\left\Vert \alpha_{k}\right\Vert _{X}^{2},\]
 using property $(iv)$ in the previous lemma in the last step. By
property $(ii)$ in the same lemma $\left\Vert \alpha_{k}\right\Vert _{X}^{2}$
is asymptotically $1,$ which shows that the first term in \ref{II pf of prop upper bound on lowe b}
tends to zero if the sequence $\mu_{k}$ is chosen as $\delta_{k}^{1/2},$
for example. To see that the second term tends to zero as well, we
estimate \[
\left\Vert (\Delta^{(k)})\alpha_{2}^{(k)}\right\Vert _{m}\leq\left\Vert (\Delta^{(k)})\alpha^{(k)}\right\Vert _{m}+\left\Vert (\Delta^{(k)})\alpha_{1}^{(k)}\right\Vert _{m}.\]
 The first term in the right hand side tends to zero by $(iii)$ in
the previous lemma and so does the second term, using \ref{II pf of prop upper bound on b b}
(that holds for any scaled sequence of forms in $\mathcal{H}_{\leq\mu_{k}}^{q}(X,L^{k}).$
This finishes the proof of the claim \ref{II pf of prop lower bound on lowe claim}.
Finally, \[
\liminf_{k}B_{\leq\mu_{k}}^{q}(0,iv/k)_{\theta_{k}}\geq\left|\alpha_{1,k}(0)\right|_{\theta_{k}}^{2}\]
and \[
\liminf_{k}\left|\alpha_{k,1,\theta_{k}}(0)\right|^{2}\geq B_{X_{0},\theta}(0,v)+0,\]
when $k$ tends to infinity, using \ref{II pf of prop lower bound on lowe claim}
and $(i)$ in the previous lemma.
\end{proof}
Now we can prove that the Morse inequalities are essentially equalities
for the space $\mathcal{H}_{\leq\mu_{k}}^{0,q}(X,L^{k})$ of low-energy
forms. But first recall that $\mathcal{H}_{\leq\mu_{k}}^{0,q}(X,L^{k})$
depends on a large parameter $R,$ since the metrics $\omega_{k}$
depend on $R.$ 

\begin{thm}
\label{thm main} Suppose that $X$ is is a compact manifold with
boundary satisfying condition $Z(q).$ Then there is a sequence $\mu_{k}$
tending to zero such that the limit of $k^{-n}\dim\mathcal{H}_{\leq\mu_{k}}^{0,q}(X,L^{k})$
when $k$ tends to infinity is equal to \[
(-1)^{q}(\frac{1}{2\pi})^{n}(\int_{X(q)}\Theta_{n}+\int_{\partial X}\int_{T(q)_{\rho,x}}(\Theta+t\mathcal{L})_{n-1}\wedge\partial\rho\wedge dt)+\epsilon_{R}\]
where the sequence $\epsilon_{R}$ tends to zero when $R$ tends to
infinity. 
\end{thm}
\begin{proof}
The proof is completely analogous to the proof of theorem \ref{II thm morse}.
In the first region one just replaces the claim \ref{II pf of first reg claim}
in the proof of proposition \ref{I prop first region} by the asymptotic
equality for $B_{\leq\mu_{k}}^{q}(0,iv/k)_{\theta_{k}}$ obtained
by combing the propositions \ref{II prop upper bound on lowenergy b}
and \ref{II prop lower bound on b low}. Moreover, a simple modification
of the proof of proposition \ref{II prop upper bound on lowenergy b}
shows that there is no contribution from the integrals over the second
and third region, when $R$ tends to infinity, as before. Finally,
the convergence on the inner part of $X$ was shown in (\cite{berm}).
\end{proof}
Recall that the Dolbeault cohomology group $H^{0,q}(X,L^{k})$ is
isomorphic to the space of harmonic forms, which is a subspace of
of $\mathcal{H}_{\leq\mu_{k}}^{0,q}(X,L^{k}).$ Hence, the previous
theorem is stronger than the weak Morse inequalities for the dimensions
$h^{q}(L^{k})$ of $H^{0,q}(X,L^{k}),$ theorem \ref{II thm morse}.
When $X$ has no boundary Demailly showed that, by combining a version
of theorem \ref{thm main} with some homological algebra, one gets
strong Morse ineqaulities for the Dolbeault cohomology groups. These
are inequalities for an alternating sum of all $h^{i}(L^{k}):$s when
the degree $i$ varies between $0$ and a fixed degree $q$ \cite{d1}.
In fact, a variation of the homological algebra argument yields inequalities
for alternating sums when the degree $i$ varies between a fixed degree
$q$ and the complex dimension $n$ of $X$ (the two versions are
related by Serre duality). However, when $X$ has a boundary one has
to impose certain curvature conditions on $\partial X$ to obtain
strong Morse inequalites from theorem \ref{thm main} . Indeed, to
apply the theorem one has to assume that $\partial X$ satisfies condition
$Z(i)$ for all degrees $i$ in the corresponding range. In particular
the corresponding dimensions will then be finite dimensional so that
the alternating sum makes sense. Now, to state the strong holomorphic
Morse inequalities for a manifold with boundary, recall that the boundary
of a compact complex manifold is called $q-$convex if the Levi form
$\mathcal{L}$ has at least $n-q$ positive eigenvalues along $T^{1,0}(\partial X)$
and it is called $q-$concave if the Levi form has at least $n-q$
negative eigenvalues along $T^{1,0}(\partial X)$ (i.e $\partial X$
is $q-$convex {}``from the inside'' precisely when it is $q-$concave
{}``from the outside''). We will denote by $X(\geq q)$ the union
of all sets $X(i)$ with $i\geq q$ and $T(\geq q)_{\rho,x}$ is defined
similarly. The sets $X(\leq q)$ and $T(\leq q)_{\rho,x}$ are defined
by puting $i\leq q$ in the previous definitions. Finally, we set
\[
I(\geq q):=(\frac{1}{2\pi})^{n}(\int_{X(\geq q)}\Theta_{n}+\int_{\partial X}\int_{T(\geq q)_{\rho,x}}(\Theta+t\mathcal{L})_{n-1}\wedge\partial\rho\wedge dt)\]
and define $I(\leq n-1-q)$ similarly. 

\begin{thm}
\label{thm:strong morse}Suppose that $X$ is an $n-$dimensional
compact manifold with boundary. If the boundary is strongly $q-$convex,
then \[
k^{-n}\sum_{i=q}^{n}(-1)^{q-j}h^{j}(L^{k})\leq I(\geq q)k^{n}+o(k^{n}).\]
If $X$ has strongly $q-$concave boundary, then \[
\sum_{i=0}^{n-1-q}(-1)^{q-j}h^{j}(L^{k})\leq I(\leq n-1-q)k^{n}+o(k^{n}).\]

\end{thm}
\begin{proof}
First note that if $\partial X$ is $q-$convex, then $\partial X$
satisfies condition $Z(i)$ for $i$ such that $n-q\leq i\leq n.$
Similarly, if $\partial X$ is $q-$concave, then $\partial X$ satisfies
condition $Z(i)$ for $i$ such that $0\leq i\leq n-q-1.$ The proof
then follows from theorem \ref{thm main} and the homological algebra
argument in \cite{d1},\cite{d2}. See also \cite{bo} and \cite{m1}. 
\end{proof}

\subsection{Strong Morse inequalities on open manifolds}

One can also define $q-$convexity and $q-$concavity on open manifolds
following Andreotti and Grauert \cite{a-g}. First, one says that
a function $\rho$ is $q-$convex if $i\partial\overline{\partial}\rho$
has at least $n-q+1$ positive eigen values. Next, an open manifold
$Y$ is said to be $q-$convex if it has an exhaustion function $\rho$
that is $q-$convex outside some compact subset $K$ of $Y.$ The
point is that the regular sublevel sets of $\rho$ are then $q-$convex
considered as compact manifolds with boundary. The extra positive
eigen value occuring in the definition of $q-$convexity for an open
manifold is needed to make sure that $i\partial\overline{\partial}\rho$
still has at least $n-q$ positive eigen values along a regular level
surface of $\rho.$ Finally, an open manifold $Y$ is said to be $q-$concave
if it has an exhaustion function $\rho$ such that $-\rho$ is $q-$convex
outside some compact subset $K$ of $Y.$ 

Now, by remark \ref{rem:open cohom} theorem \ref{thm:strong morse}
extends to any $q-$convex open manifold $Y$ with a line bundle $L$
if one uses the usual Dolbeault cohomology $H^{0,*}(Y,L^{k})$ (or
equivalently the sheaf cohomology $H^{*}(Y,\mathfrak{\mathcal{O}}(L^{k})$)
and the curvature integrals are taken over a regular level surface
of $\rho$ in the complement of the compact set $K.$ However, for
a $q-$concave open manifold $Y$ one only gets the corresponding
result if $n-q-1$ is replaced with $n-q-2.$ Indeed, by remark \ref{rem:open cohom}
one has to make sure that condition $Z(i+1)$ holds for the highest
degree $i$ occuring in the alternating sum. In this form the $q-$convex
case and $q-$concave case was obtained by Bouche \cite{bo} and Marinescu
\cite{m1}, respectively, under the assumption that the curvature
of the line bundle $L$ is adapted to the curvature of the boundary
of $X$ in a certain way. Comparing with theorem \ref{thm:strong morse}
their assumptions imply that the boundary integral vanishes. There
is also a very recent preprint \cite{m2} of Marinescu where strong
Morse inequalities on a $q-$concave manifold with an arbitrary line
bundle $L$ are obtained. However, the corresponding boundary term
is not as precise as the one in theorem \ref{thm:strong morse} and
in section \ref{sec:Sharp} we will show that theorem \ref{thm:strong morse}
is sharp.

Note that since the curvature integrals are taken over \emph{any}
regular level surface of $\rho$ in $Y$ one expects that $I(\geq q)$
and $I(\leq n-1-q)$ are independent of the level surface. This is
indeed the case (see remark \ref{rem:indep of right hand side}).

\subsection{Application to the volume of semi-positive line bundles}

The most interesting case when the strong Morse inequalities above
apply is when $X$ is a strongly pseudoconcave manifold $X$ of dimension
$n\geq3$ with a semi-positive line bundle $L,$ i.e the Levi form
$\mathcal{L}$ is negative along along $T^{1,0}(\partial X)$ and
the curvature form of $L$ is semi-positive in $X.$ Then one gets
a lower bound on the dimension of the space of holomorphic sections
with values in $L^{k}.$ Namely, $h^{0}(L^{k})$ is asymptotically
bounded from below by\begin{equation}
(\frac{1}{2\pi})^{n}(\int_{X(0)}\Theta_{n}+\int_{\partial X}\int_{T(\leq1)}(\Theta+t\mathcal{L})_{n-1}\wedge\partial\rho\wedge dt)k^{n}+h^{1}(L^{k})+o(k^{n})\label{eq strong morse}\end{equation}
 In particular, if the curvatures are such that the coefficient in
front of $k^{n}$ is positive, then the dimension of $H^{0}(X,L^{k})$
grows as $k^{n}.$ In other words, the line bundle $L$ is \emph{big}
then. For example, this happens when the curvature forms are conformally
equivalent along the complex tangential directions, i.e. if there
is a function $f$ on $\partial X$ such that \begin{equation}
\mathcal{L}=-f\Theta\label{eq conformal reflation}\end{equation}
when restricted to $T^{1,0}(\partial X)\otimes T^{0,1}(\partial X).$
In fact, by multiplying the original $\rho$ by $f^{-1}$ we may and
will assume that $f=1.$ The lower bound \ref{eq strong morse} combined
with the upper bound from the weak Morse inequalities (theorem \ref{II thm morse})
then gives the following corollary.

\begin{cor}
Suppose that $X$ is a strongly pseudoconcave manifold $X$ of dimension
$n\geq3$ with a semi-positive line bundle $L.$ Then if the curvature
forms are conformally equivalent at the boundary \[
h^{0}(L^{k})=k^{n}(\frac{1}{2\pi})^{n}(\int_{X}\Theta_{n}+\frac{1}{n}\int_{\partial X}(i\partial\overline{\partial}\rho)_{n-1}\wedge i\partial\rho)+o(k^{n}).\]

\end{cor}
When $L$ is positive, the conformal equivalence in the previous corollary
says that the symplectic structure on $X$ determined by $L$ is compatible
with the contact structure of $\partial X$ determined by the complex
structure, in a strong sense (compare \cite{e}) and the conclusion
of the corollary may be expressed by the formula\begin{equation}
\textrm{Vol$(L)=$Vol$(X)+\frac{1}{n}$Vol$(\partial X)$}\label{II eq vol is vol}\end{equation}
in terms of the symplectic and contact volume of $X$ and $\partial X,$
respectively (where the volume of a line bundle $L$ is defined as
the lim sup of $(2\pi)^{n}k^{-n}h^{0}(L^{k})).$ The factor $\frac{1}{n}$
in the formula is related to the fact that if $(X_{+},d\alpha)$ is
a $2n-$dimensional real symplectic manifold with boundary, such that
$\alpha$ is a contact form for $\partial X_{+},$ then, by Stokes
theorem, the contact volume of $\partial X_{+}$ divided by $n$ is
equal to the symplectic volume of $X_{+}.$ In fact, this is how we
will show that \ref{II eq vol is vol} is compatible with hole filling
in the following section.

\section{\label{sec:Sharp}Sharp examples and hole filling}

In this section we will show that the leading constant in the Morse
inequalities \ref{II thm morse} is sharp. When $X$ is a compact
manifold without boundary, this is well-known. Indeed, let $X$ be
the $n-$dimensional flat complex torus $\C^{n}/\Z^{n}+i\Z^{n}$ and
consider the hermitian holomorphic line bundle $L_{\lambda}$ over
$T^{n}$ determined by the constant curvature form \[
\Theta=\sum_{i=1}^{n}\frac{i}{2}\lambda_{i}dz_{i}\wedge\overline{dz_{i}},\]
 where $\lambda_{i}$ are given non-zero integers \cite{gri}. Then
one can show (see the remark at the end of the section) that \begin{equation}
B^{q}(x)\equiv\frac{1}{\pi^{n}}1_{X(q)}\left|\textrm{det}_{\omega}\Theta\right|,\label{II example compact b}\end{equation}
 where $1_{X(q)}$ is identically equal to one if exactly $q$ of
the eigenvalues $\lambda_{i}$ are negative and equal to zero otherwise.
This shows that the leading constant in the Morse inequalities on
a compact manifold is sharp. 

Let us now return to the case of a manifold with boundary. We let
$X$ be the manifold obtained as the total space of the unit disc
bundle in the dual of the line bundle $L_{\mu}$ (where $L_{\mu}$
is defined as above) over the torus $T^{n-1},$ where $\mu_{i}$ are
$n-1$ given non-zero integers. Next, we define a hermitian holomorphic
line bundle over $X.$ Denote by $\pi$ the natural projection from
$X$ onto the torus $T^{n-1}.$ Then the pulled back line bundle $\pi^{*}L_{\lambda}$
is a line bundle over $X.$ The construction is summarized by the
following commuting diagram \[
\begin{array}{ccccc}
\pi^{*}L_{\lambda} &  &  &  & L_{\lambda}\\
\downarrow &  &  &  & \downarrow\\
X & \hookrightarrow & L_{\mu}^{*} & \rightarrow & T^{n-1}\end{array}\]
Let $h$ be the positive real-valued function on $X,$ defined as
the restriction to $X$ of the squared fiber norm on $L_{\mu}^{*}.$
Then $\rho:=\ln h$ is a defining function for $X$ close to the boundary
and we define a hermitian metric $\omega$ on $X$ by \[
\omega=\frac{i}{2}\partial\overline{\partial}\left|z\right|^{2}+\frac{i}{2}h^{-1}\partial h\wedge\overline{\partial h}\]
 extended smoothly to the base $T^{n-1}$ of $X.$ The following local
description of the situation is useful. The part of $X$ that lies
over a fundamental domain of $T^{n-1}$ can be represented in local
holomorphic coordinates $(z,w),$ where $w$ is the fiber coordinate,
as the set of all $(z,w)$ such that \[
h(z,w)=\left|w\right|^{2}\exp(+\sum_{i=1}^{n-1}\mu_{i}\left|z_{i}\right|^{2})\leq1\]
 and the fiber metric $\phi$ for the line bundle $\pi^{*}L_{\lambda}$
over $X$ may be written as \[
\phi(z,w)=\sum_{i=1}^{n-1}\lambda_{i}\left|z_{i}\right|^{2}.\]
The proof of the following proposition is very similar to the proof
of theorem \ref{II thm model}, but instead of Fourier transforms
we will use Fourier series, since the $\R-$ symmetry is replaced
by an $S^{1}-$ symmetry (the model domain $X_{0}$ in section \ref{II section model}
is the universal cover of $X$ defined above). 

\begin{thm}
\label{thm: model}Let $J(q)$ be the set of all integers $j$ such
that the form $\partial\overline{\partial}\phi+j\partial\overline{\partial}\rho$
has exactly $q$ negative eigenvalues. Then \begin{equation}
B_{X}^{q}=(\frac{1}{2\pi})^{n}\sum_{j\in J(q)}\textrm{det}_{\omega}(i(\partial\overline{\partial}\phi+j\partial\overline{\partial}\rho))\frac{1}{2}(j+1)h^{j}.\label{II prop example i}\end{equation}
In particular, the dimension of $H^{0,q}(X,\pi^{*}L_{\lambda})$ is
given by \begin{equation}
(\frac{i}{2\pi})^{n}\int_{\partial X}\sum_{j\in J(q)}(\partial\overline{\partial}\phi+j\partial\overline{\partial}\rho)_{n-1}\wedge\partial\rho\label{II prop example ii }\end{equation}
 and the limit of the dimensions of $H^{0,q}(X,(\pi^{*}L_{\lambda})^{k})$
divided by $k^{n}$ is \begin{equation}
(\frac{i}{2\pi})^{n}\int_{\partial X}\int_{T_{x,\rho}(q)}(\partial\overline{\partial}\phi+t\partial\overline{\partial}\rho)_{n-1}\wedge\partial\rho\wedge dt\label{II prop ex iii}\end{equation}

\end{thm}
\begin{proof}
First note that if $\alpha_{j}$ is a form on $T^{n-1}$ with values
in $L_{\mu}^{j}\otimes L_{\lambda},$ then \[
\alpha(z,w):=\sum_{j\geq0}\alpha_{j}(z)w^{j}\]
 defines a global form on $X$ with values in $\pi^{*}L_{\lambda}.$
The proof of proposition \ref{II prop tangential} can be adapted
to the present situation to show that any form $\alpha$ in $\mathcal{H}^{0,q}(X,\phi)$
is of this form with $\alpha_{j}$ in $\mathcal{H}^{0,q}(T^{n-1},L_{\mu}^{j}\otimes L_{\lambda}).$
Actually, since $X$ is a fiber bundle over $T^{n-1}$ with \emph{compact}
fibers one can also give a somewhat simpler proof. For example, to
show that $\alpha$ is tangential one solves the $\overline{\partial}-$
equation along the fibers of closed discs in order to replace the
normal part $\alpha_{N}$ with an exact form. Then using the assumption
the $\alpha$ is coclosed one shows that the exact form must vanish.
The details are omitted. We now have the following analog of proposition
\ref{II lemma model scalar pr} for any $\alpha$ in $\mathcal{H}^{0,q}(X,\phi)$
\begin{equation}
(\alpha,\alpha)=2\pi\sum_{j}(\alpha_{j},\alpha_{j})b_{j},\,\,\, b_{j}=\int_{0}^{1}(r^{2})^{j}rdr=1/2(j+1)^{-1}\label{II sharp ex sclar pr}\end{equation}
 in terms of the induced norms. To see this, one proceed as in the
proof of proposition \ref{II lemma model scalar pr}, now using the
Taylor expansion of $\alpha.$ Writing $\psi(z)=\sum_{i=1}^{n-1}\mu_{i}\left|z_{i}\right|^{2}$
and restricting $z$ to the fundamental region of $T^{n-1}$ we get
that $(\alpha,\alpha)$ is given by \[
\int_{\left|w\right|^{2}<e^{-\psi(z)}}\left|\sum_{j}\alpha_{j}(z)w^{j}\right|^{2}e^{-\phi(z)}(\frac{i}{2}\partial\overline{\partial}\left|z\right|^{2})_{n-1}e^{\psi(z)}rdrd\theta.\]
 Now using Parceval's formula for Fourier series in the integration
over $\theta$ this can be written as \[
2\pi\sum_{j}\int_{z}\left|\widehat{\alpha_{j}}(z)\right|^{2}e^{-\phi(z)}(\frac{i}{2}\partial\overline{\partial}\left|z\right|^{2})_{n-1}e^{\psi(z)}\int_{0}^{e^{-\psi(z)/2}}(r^{2})^{j}rdr.\]
 Finally, the change of variables $r'=e^{\psi(z)/2}r$ in the integral
over $r$ gives a factor $e^{-j\psi(z)}$ and the upper integration
limit becomes $1.$ This proves \ref{II sharp ex sclar pr}. 

As in the proof of theorem \ref{II thm model} we infer that $B_{X}$
may be expanded as \[
B_{X}(z,w)=\frac{1}{2\pi}\sum_{j}B_{j}(z)h^{j}b_{j}^{-1},\]
 where $B_{j}$ is the Bergman function of the space $\mathcal{H}^{0,q}(T^{n-1},L_{\mu}^{j}\otimes L_{\lambda}).$
According to \ref{II example compact b}, we have that \[
B_{j}(z)\equiv(\frac{1}{2\pi})^{n-1}\delta_{j,J(q)}\left|\textrm{det}_{\omega}(i(\partial\overline{\partial}\phi+j\partial\overline{\partial}\rho))\right|,\]
 where $J(q)=\{ j:\,\textrm{index$($}\partial\overline{\partial}\phi+j\partial\overline{\partial}\rho)=q\}$and
where the sequence $\delta_{j,J(q)}$ is equal to $1$ if $j\in J(q)$
and zero otherwise. Thus, \ref{II prop example i} is obtained. Integrating
\ref{II prop example i} over $X$ gives \[
\int_{X}B_{X}\omega_{n}=\frac{1}{2\pi}\int_{T^{n-1}}\sum_{j}B_{j}(\frac{i}{2}\partial\overline{\partial}\left|z\right|^{2})_{n-1}\int_{0}^{2\pi}d\theta\int_{0}^{1}(r^{2})^{j}rdrb_{j}^{-1}.\]
 The integral over the radial coordinate $r$ is cancelled by $b_{j}^{-1}$
and we may write the resulting integral as \[
\frac{1}{2\pi}\int_{\partial X}\sum_{j}B_{j}\omega_{n-1}\wedge i\partial\rho=(\frac{i}{2\pi})^{n}\int_{\partial X}\sum_{j\in J(q)}(\partial\overline{\partial}\phi+j\partial\overline{\partial}\rho)_{n-1}\wedge\partial\rho\]
 Hence, \ref{II prop example ii } is obtained. Finally, applying
the formula \ref{II prop example ii } to the line bundle $(\pi^{*}L_{\lambda})^{k}=\pi^{*}(L_{\lambda}^{k})$
shows, since the curvature form of $\pi^{*}(L_{\lambda}^{k})$ is
equal to $k\partial\overline{\partial}\phi,$ that \[
k^{-n}\int_{X}B_{X}^{q,k}\omega_{n}=(\frac{i}{2\pi})^{n}\int_{\partial X}\sum_{j}(\partial\overline{\partial}\phi+\frac{j}{k}\partial\overline{\partial}\rho)^{n-1}\frac{1}{k}\wedge\partial\rho,\]
 where the sum is over all integers $j$ such that $\partial\overline{\partial}\phi+\frac{j}{k}\partial\overline{\partial}\rho$
has exactly $q$ negative eigenvalues. Observe that the sum is a Riemann
sum and when $k$ tends to infinity we obtain \ref{II prop ex iii}. 
\end{proof}
Note that since the line bundle $\pi^{*}L_{\lambda}$ over $X$ is
flat in the fiber direction the integral over $X$ in \ref{main statement}
vanishes. Hence, the theorem above shows that the holomorphic Morse
inequalities are sharp. The most interesting case covered by the theorem
above is when the line bundle $\pi^{*}L_{\lambda}$ (simply denoted
by $L)$ over $X$ is semi-positive, and positive along the tangential
directions, and $X$ is strongly pseudoconcave. This happens precisely
when all $\lambda_{i}$ are positive and all $\mu_{i}$ are negative.
Then, for $n\geq3,$ the theorem above shows that the dimension of
$H^{0,1}(X,L^{k})$ grows as $k^{n}$ unless the curvature of $L$
is a multiple of the Levi curvature of the boundary, i.e. unless $\lambda$
and $\mu$ are parallel as vectors. This is in contrast to the case
of a manifold without boundary, where the corresponding growth is
of the order $o(k^{n})$ for a semi-positive line bundle. Note that
the bundle $L$ above always admits a metric of \emph{positive} curvature.
Indeed, the fiber metric $\phi+\epsilon h$ on $L$ can be seen to
have positive curvature, if the positive number $\epsilon$ is taken
sufficiently small. However, if $\lambda$ and $\mu$ are not parallel
as vectors, there is no metric of positive curvature which is conformally
equivalent to the Levi curvature at the boundary. This follows from
the weak holomorphic Morse inequalities, theorem \ref{II thm morse},
since the growth of the dimensions of $H^{0,1}(X,L^{k})$ would be
of the order $o(k^{n})$ then.

\begin{rem}
To get examples of open manifold $Y$ as described in remark \ref{rem:open cohom}
one may take the total space of the line bundle $L_{\mu}^{*}$ over
$T^{n-1},$ as defined in the beginning of the section. Then $\rho$
is an exhaustion sequence, exhausting $L_{\mu}^{*}$ by disc bundles.
Furthermore, to get examples of manifolds with boundary $X$ where
the index of the Levi curvaure form is non-constant one may take $X$
to be an annulus bundle in $L_{\mu}^{*}.$ Such a manifold is neither
$q-$convex or $q-$concave for any $q.$ Theorem \ref{thm: model}
extends to such manifolds $X$ if one uses Laurent expansions of sections
instead of Taylor expansions. A concrete example is given by the hyperplane
bundle $O(1)$ over $\P^{n-1}.$ Then the corresponding annulus bundle
is biholomorphic to a spherical shell in $\C^{n-1},$ i.e. all $z$
in $\C^{n-1}$ such that $r\leq\left|z\right|\leq r'$ for some given
numbers $r$ and $r'.$ It has one pseudoconvex and one pseudoconcave
boundary component.
\end{rem}
Finally, a remark about the proof of formula \ref{II example compact b}.

\begin{rem}
To prove formula \ref{II example compact b} one can for example reduce
the problem to holomorphic sections, i.e. when $q=0$ (compare \cite{berm}).
One could also use symmetry to first show that the Bergman kernel
is constant and then compute the dimension of $H^{q}(T^{n},L_{\lambda})$
by standard methods. To compute the dimension one writes the line
bundle $L_{\lambda}$ as $L_{\lambda}=\pi_{1}^{*}L^{\lambda_{1}}\otimes\pi_{2}^{*}L^{\lambda_{2}}\otimes...\otimes\pi_{n}^{*}L^{\lambda_{n}},$
using projections on the factors of $T^{n},$ where $L$ is the classical
line bundle over the elliptic curve $T^{1}=\C/\Z+i\Z$, such that
$H^{0}(\C/\Z+i\Z,L)$ is generated by the Riemann theta function \cite{gri}.
Now, using Kunneth's theorem one gets that $H^{q}(T^{n},L_{\lambda})$
is isomorphic to the direct sum of all tensor products of the form\[
H^{1}(T^{1},L^{\lambda_{i_{1}}})\otimes\cdots\otimes H^{1}(T^{1},L^{\lambda_{i_{q}}})\otimes H^{0}(T^{1},L^{\lambda_{i_{q+1}}})\otimes\cdots\otimes H^{0}(T^{1},L^{\lambda_{i_{n}}}).\]
Observe that this product vanishes unless the index $I=(i_{1},...,i_{n})$
is such that the first $q$ indices are negative while the others
are positive. Indeed, first observe that if $m$ is a positive integer,
the dimension of $H^{0}(T^{1},L^{-m})$ vanishes, since $L^{-m}$
is a negative line bundle. Next, by Serre duality $H^{1}(T^{1},L^{m})\approxeq H^{0}(T^{1},L^{-m}),$
since the canonical line bundle on $T^{1}$ is trivial. So the dimension
of $H^{1}(T^{1},L^{m})$ vanishes as well. In particular, the dimension
of $H^{q}(T^{n},L_{\lambda})$ vanishes unless exactly $q$ of the
numbers $\lambda_{i}$ are negative, i.e. unless the index of the
curvature of $L_{\lambda}$ is equal to $q.$ Finally, if the index
is equal to $q,$ then, using that $H^{0}(T^{1},L)$ is one-dimensional,
combined with Serre duality and Kunneth's formula again, one gets
that the dimension of $H^{0}(T^{1},L^{-1})$ is equal to the absolute
value of the product of all eigenvalues $\lambda_{i}.$ This proves
\ref{II example compact b}.
\end{rem}

\subsection{\label{section hole and contact}Relation to hole filling and contact
geometry. }

Consider a compact strongly pseudoconcave manifold $X$ with a semi-positive
line bundle $L.$ We will say that the pair $(X,L)$ may be \emph{filled}
if there is a compact complex manifold $\widetilde{X},$ without boundary,
with a semi-positive line bundle $\widetilde{L}$ such that there
is a holomorphic line bundle injection of $L$ into $\widetilde{L}.$
\footnote{By a theorem of Rossi \cite{r}, the pair $(X,L)$ may always be filled
if $L$ is the trivial holomorphic line bundle.%
} The simplest situation is as follows. Start with a compact complex
manifold $\widetilde{X}$ with a positive line bundle $\widetilde{L}$
(by the Kodaira embedding theorem $\widetilde{X}$ is then automatically
a projective variety \cite{gri}). We will obtain a pseudoconcave
manifold $X$ by making a small hole in $\widetilde{X}.$ Consider
a small neighborhood of a fixed point $x$ in $\widetilde{X},$ holomorphically
equivalent to a ball in $\C^{n},$ where $\widetilde{L}$ is holomorphically
trivial and let $\phi$ be the local fiber metric. We may assume that
$\phi(x)=0$ and that $\phi$ is non-negative close to $x.$ Then
for a sufficiently small $\epsilon$ the set where $\phi$ is strictly
less than $\epsilon$ is a strongly pseudoconvex domain of $\widetilde{X}$
and its complement is then a strongly pseudoconcave manifold that
we take to be our manifold $X.$ We let $L$ be the restriction of
$\widetilde{L}$ to $X.$ A defining function of the boundary of $X$
can be obtained as $\rho=-\phi.$ Now, since $\widetilde{L}$ is a
positive line bundle it is wegv thesischarles

ll-known that \[
\lim_{k}k^{-n}\dim_{\C}H^{0}(\widetilde{X},\widetilde{L}^{k})=(\frac{i}{2\pi})^{n}\frac{1}{n!}\int_{\widetilde{X}}(\partial\overline{\partial}\widetilde{\phi})^{n}.\]
 In fact, this holds for any semi-positive line bundle, as can be
seen by combining Demailly's holomorphic Morse inequalities \ref{demailly weak ineaq}
with the Riemann-Roch theorem (this was first proved by different
methods in \cite{s2}).

On the other hand we have, by Harthog's phenomena, that $H^{0}(\widetilde{X},\widetilde{L}^{k})$
is isomorphic to $H^{0}(X,L^{k}).$ So decomposing the integral above
with respect to \begin{equation}
\widetilde{X}=X\bigsqcup X^{c}\label{II example hole decomp}\end{equation}
 and using Stokes theorem gives that \begin{equation}
\lim_{k}k^{-n}\dim_{\C}H^{0}(X,L^{k})=(\frac{i}{2\pi})^{n}(\frac{1}{n!}\int_{X}(\partial\overline{\partial}\phi)_{n}-\frac{1}{n!}\int_{\partial X}(\partial\overline{\partial}\phi)^{n-1}\wedge\partial\phi)\label{II ex hole dim}\end{equation}
 Let us now compare the boundary integral above with the curvature
integral in the holomorphic Morse inequalities \ref{main statement}.
Since $\rho=-\phi$ this integral equals \[
-\frac{1}{(n-1)!}\int_{\partial X\times[0,1]}((1-t)\partial\overline{\partial}\phi)^{n-1}\wedge\partial\phi\wedge dt,\]
 which coincides with the boundary integral in \ref{II ex hole dim}
since $\int_{o}^{1}(1-t)^{n-1}dt=1/n.$ This shows that the holomorphic
Morse inequalities, theorem \ref{thm main} are sharp for the line
bundle $L$ over $X.$ To show that the Morse inequalities are sharp
as soon as a pair $(X,L)$ may be filled by a Stein manifold it is
useful to reformulate the boundary term in \ref{main statement} in
terms of the contact geometry of the boundary $\partial X.$

Let us first recall some basic notions of contact geometry (\cite{ar}).
The distribution $T^{1,0}(\partial X)$ can be obtained as $\ker(-i\partial\rho)$
and since, by assumption, the restriction of $d(-i\partial\rho)$
is non-degenerate it defines a so called \emph{contact} distribution
and $\partial X$ is hence called a \emph{contact manifold}. By duality
$T^{1,0}(\partial X)$ determines a real line bundle in the real cotangent
bundle $T^{*}(\partial X)$ that can be globally trivialized by the
form $-i\partial\rho.$ Denote by $X_{+}$ the associated fiber bundle
over $\partial X$ of {}``positive'' rays and denote by $\alpha$
the tautological one form on $T^{*}(\partial X),$ so that $d\alpha$
is the standard symplectic form on $T^{*}(\partial X).$ The pair
$(X_{+},d\alpha)$ is called the \emph{symplectification} of the contact
manifold $\partial X$ in the literature \cite{ar}. More concretely,
\[
X_{+}=\left\{ t(-i\partial\rho_{x}):\, x\in\partial X,t\geq0\right\} ,\]
i.e. $X_{+}$is isomorphic to $\partial X\times[0,\infty[$ and $\alpha=-it\partial\rho$
so that $d\alpha=i(t\partial\overline{\partial}\rho+\partial\rho\wedge dt).$
The boundary integral in \ref{main statement} may now be compactly
written as \[
\int_{X_{+}(q)}(\Theta+d\alpha)_{n},\]
 where $X_{+}(q)$ denotes the part of $X_{+}$ where the pushdown
of $d\alpha$ to $\partial X$ has exactly $q$ negative eigenvalues
along the contact distribution $T^{1,0}(\partial X).$ 

Let us now assume that $X$ is strongly pseudoconcave and that $(X,L)$
is filled by $(\widetilde{X},L)$ (abusing notation slightly). We
will also assume that the strongly pseudoconvex manifold $Y,$ in
$\widetilde{X},$ obtained as the closure of the complement of $X$
in $\widetilde{X},$ has a defining function that we write as $-\rho$
which is plurisubharmonic on $Y.$ We may assume that the set of critical
points of $-\rho$ on $Y$ is finite and to simplify the argument
we assume that there is exactly one critical point $x_{0}$ in $Y$
and we assume that $\rho(x_{0})=1$ (the general argument is similar).
For a regular value $c$ of $\rho$ we let $X_{+}(0)_{c}$ be the
subset of the symplectification of $\rho^{-1}(c)$ defined as above,
thinking of $\rho^{-1}(c)$ as a strictly pseudoconcave boundary.
Now consider the following manifold with boundary: \[
\mathcal{X}_{\varepsilon}(0)=\bigcup_{c\in[0,1-\varepsilon]}X_{+}(0)_{c}.\]
 More concretely, $\mathcal{X}_{\varepsilon}(0)$ can be identified
with a subset of the positive closed cone in $T^{*}(Y,\C)$ determined
by $\partial\rho:$ \[
\left\{ t(-i\partial\rho_{x}):\, x\in Y,t\geq0\right\} .\]
 Hence, $\mathcal{X}_{\varepsilon}(0)$ is a fiber bundle over a subset
of $Y$ and when $\epsilon$ tends to zero, the base of $\mathcal{X}_{\varepsilon}(0)$
tends to $Y.$ Note that the fibers of $\mathcal{X}_{\varepsilon}(0)$
are a finite number of intervals and the induced function $t$ on
$\mathcal{X}_{\varepsilon}$ is uniformly bounded with respect to
$\varepsilon$ (i.e. the {}``height'' of the fiber is uniformly
bounded). Indeed, we have assumed that $i\partial\overline{\partial}\rho$
is strictly negative. This forces $\Theta+ti\partial\overline{\partial}\rho$
to be negative on all of $Y$ for all $t$ larger then some fixed
number $t_{0}.$ In particular such a $t$ is not in $T_{x}(0)$ for
any $x$ in $Y,$ i.e. not in any fiber of $\mathcal{X}_{\varepsilon}(0).$
Now observe that the form $\Theta+d\alpha$ on $X_{+}(0)$ extends
to a closed form in $\mathcal{X}_{\varepsilon}(0)$ and the restriction
of the form to $Y$ coincides with $\Theta.$ Let us now integrate
the form $(\Theta+d\alpha)_{n}$ over the boundary of $\mathcal{X}_{\varepsilon}(0).$
The boundary can be written as \[
\partial(\mathcal{X}_{\varepsilon}(0))=X_{+}(0)\bigcup\left(\bigcup_{c\in]0,1-\varepsilon[}\partial(X_{+}(0)_{c})\right)\bigcup X_{+}(0)_{\varepsilon}.\]
Since the form is closed, the integral over $\partial(\mathcal{X}_{\varepsilon}(0))$
vanishes according to Stokes theorem, giving \[
0=\int_{X_{+}(0)}(\Theta+d\alpha)_{n}-\int_{Y}\Theta{}_{n}+0+O(\varepsilon)\]
 where the zero contribution comes from the fact that the form $(\Theta+d\alpha)_{n}$
vanishes along $\left(\bigcup_{c\in]0,1-\varepsilon[}\partial(X_{+}(0)_{c})\right)-Y.$
The term $O(\varepsilon)$ comes from the integral (of a uniformly
bounded function) over the {}``cylinder'' $X_{+}(0)_{\varepsilon}$
around the point $x_{0}.$ Finally, by letting $\epsilon$ tend to
zero we see that the Morse inequalities for $L$ over $X$ are sharp
in this situation as well.

\begin{rem}
\label{rem:indep of right hand side}The preceeding argument also
shows that if $\rho$ is a function on an open manifold $Y$ with
regular values $c$ and $c'$ (where $c$ is less than $c'),$ then
\[
\int_{X_{+}(i)_{c}}(\Theta+d\alpha)_{n}=\int_{\rho^{-1}]c,c']}\Theta_{n}+\int_{X_{+}(i)_{c'}}(\Theta+d\alpha)_{n}\]
 for all $i$ such that $i\geq q,$ if $\rho$ is $q-$convex on $\rho^{-1}]c,c'].$
In other words, the right hand sides in the weak Morse inequalities
for $\rho^{-1}(\leq c)$ and $\rho^{-1}(\leq c')$ coincide. The analogous
statement also holds in the $q-$convave case.
\end{rem}
\begin{acknowledgement}
The author is grateful to his advisor Bo Berndtsson for for many enlightening
discussions and fruitful suggestions during the preparation of the
paper.
\end{acknowledgement}


\begin{thebibliography}{10}
\bibitem{a}Andreotti, A: Theoremes de dependence algebrique sur les espaces complexes
pseudoconcaves. Bull. Soc. Math. France, 91, 1963, 1--38. 
\bibitem{a-g}Andreotti, A, Grauert, H: Theoremes de finitude pour la cohomologie
des espaces complexes, Bull. Soc. Math. France, 90, 1962, 193--259. 
\bibitem{ar}Arnold, V.I: Symplectic Geometry. In Dynamical systems IV, 1--138
Encyclopaedia Math. Sci, 4, Springer, Berlin, 2001
\bibitem{berm}Berman, R: Bergman kernels and local holomorphic Morse inequalities.
Math Z., Vol 248, Nr 2 (2004), 325--344 (arXiv.org/abs/math.CV/0211235) 
\bibitem{berm2}Berman, R: Super Toeplitz operators on holomorphic line bundles (arXiv.org/abs/math.CV/0406032) 
\bibitem{bern}Berndtsson, Bo: Bergman kernels related to hermitian line bundles
over compact comlex manifolds. Explorations in complex and Riemannian
geometry, 1--17, Contemp. Math, 332, Amer. Math. Soc, Providence,
RI, 2003. 
\bibitem{bo}Bouche, T: Inegalite de Morse pour la d''-cohomologie sur une variete
non-compacte. Ann. Sci. Ecole Norm. Sup. 22, 1989, 501-513
\bibitem{cm}Chern, S.S; Moser.J.K: Real hypersurfaces in complex manifolds, Acta
Math. 133 (1974), 219--271
\bibitem{d1}Demailly, J-P: Champs magnetiques et inegalite de Morse pour la d''-cohomologie.,
Ann Inst Fourier, 355 (1985,185-229)
\bibitem{d2}Demailly, J-P: Holomorphic Morse inequalities. Several complex variables
and complex geometry, Part 2 (Santa Cruz, CA, 1989), 93-114
\bibitem{d3}Demailly, J-P: Introduction à la théorie de Hodge. In {}``Transcendental
methods in algebraic geometry. Lectures given at the 3rd C.I.M.E.
Session held in Cetraro, July 4--12, 1994.'' Lecture Notes in Mathematics,
1646. Springer-Verlag, 1996.
\bibitem{10}Donnelly, H, Fefferman, $L^{2}-$ cohomology and index theorem for
the Bergman metric. Ann. of Math. (2) 118, 1983, no 3. 593-618
\bibitem{e}Eliashberg, Y: A few remarks about symplectic filling. Geometry and
topology, Vol 8 (2004) Nr 6, 277---293
\bibitem{f-k}Folland, G.B, Kohn J.J:The Neumann problem for the Cauchy-Riemann
complex. Annals of Math. Studies 75, Princeton University Press, 1972. 
\bibitem{gri}Griffiths, P; Harris, J: Principles of algebraic geometry. Wiley Classics
Library. John Wiley \& Sons, Inc., New York, 1994.
\bibitem{gro}Gromov, M: Kähler hyperbolicity and $L^{2}$-Hodge theory. J. Differential
Geom. 33 (1991), no. 1, 263--292. 
\bibitem{h}Hörmander, L: $L^{2}$ estimates and existence theorems for the $\overline{\partial}-$operator.
Acta Math. 113 1965 89--152.
\bibitem{m1}Marinescu, G: Asymptotic Morse inequalities for Pseudoconcave manifolds.
Ann. Scuola. Norm. Sup. Pisa CL Sci (4) 23 (1996), no 1, 27--55
\bibitem{m2}Marinescu, G: Existence of holomorphic sections and perturbation of
positive line bundles over $q-$concave manifolds. (arXiv.org/abs/math.CV/0402041) 
\bibitem{r}Rossi, H: Attaching analytic spaces to an analytic space along a pseudoconcave
boundary. Proc. Conf. Complex. Manifolds (Minneapolis), Springer-Verlag,
New York 1965, 242--256
\bibitem{ru}Rudin, W: Real and complex analysis. McGraw-Hill Book Company, international
edition 1987.
\bibitem{s}Siu, Y. T: Some recent results in complex manifold theory related
to vanishing theorems for the semipositive case. Workshop Bonn 1984
(Bonn, 1984), 169--192, Lecture Notes in Math., 1111, Springer, Berlin,
1985.
\bibitem{s2}Siu, Y.T: A vanishing theorem for semipositive line bundles over non-Kähler
manifolds. J. Differential Geom. 19 (1984), no. 2, 431--452. 
\bibitem{we}Wells, R. O., Jr.: Differential analysis on complex manifolds. Graduate
Texts in Mathematics, 65. Springer-Verlag, New York-Berlin, 1980. 
\bibitem{wi}Witten, E: Supersymmetry and Morse theory. J. Differential Geom. 17
(1982), no. 4, 661--692. \end{thebibliography}
\end{document}